\pdfoutput=1

\documentclass[12pt]{article}

\usepackage{amsmath}
\usepackage{amsthm}
\usepackage{amssymb}
\usepackage[noadjust]{cite}
\usepackage{url}
\usepackage[multiple]{footmisc}
\usepackage{graphicx}
\usepackage{epstopdf}
\usepackage{chngcntr}
\usepackage{enumitem}
\usepackage{hhline}
\usepackage{array}
\usepackage{hyperref}
\usepackage{color}
\usepackage{array}
\usepackage{multirow}
\usepackage{amsmath}
\usepackage{mathtools}
\usepackage{makecell}
\usepackage{booktabs}
\usepackage{siunitx}

\usepackage{cleveref}
\usepackage{aliascnt} 

\usepackage{mdframed}

\usepackage{color,soul}

\usepackage{tikz}

\usepackage{xcolor}

\usepackage[ruled,linesnumbered]{algorithm2e}
\SetKwInput{KwData}{\textbf{Input}}
\SetKwInput{KwResult}{\textbf{Output}}

\let\oldnl\nl
\newcommand{\nonl}{\renewcommand{\nl}{\let\nl\oldnl}}

\usepackage{float}
\floatname{algorithm}{Algorithm}

\hypersetup{hidelinks}

\makeatletter
\let\@fnsymbol\@arabic
\makeatother

\if@twoside
\else 
\fi

\numberwithin{equation}{section}
\counterwithin{figure}{section}
\counterwithin{table}{section}

\makeatletter
\let\c@table\c@figure
\makeatother

\newcommand{\newsharedthm}[2]{%
  \newaliascnt{#1}{theorem}%
  \newtheorem{#1}[#1]{#2}%
  \aliascntresetthe{#1}%
}
\theoremstyle{plain}
\newtheorem{theorem}{Theorem}[section]
\newsharedthm{lemma}{Lemma}
\newsharedthm{proposition}{Proposition}
\newsharedthm{corollary}{Corollary}
\newtheorem{problem}{Problem}

\theoremstyle{definition}
\newsharedthm{definition}{Definition}
\newsharedthm{assumption}{Assumption}
\newsharedthm{example}{Example}
\newsharedthm{remark}{Remark}
\newtheorem*{notation}{Notation}

\theoremstyle{plain}
\newtheorem{probleminner}{Problem}

\crefname{theorem}{Theorem}{Theorems}
\crefname{definition}{Definition}{Definitions}
\crefname{assumption}{Assumption}{Assumptions}
\crefname{example}{Example}{Examples}
\crefname{lemma}{Lemma}{Lemmas}
\crefname{proposition}{Proposition}{Propositions}
\crefname{corollary}{Corollary}{Corollaries}
\crefname{figure}{Figure}{Figures}
\crefname{section}{Section}{Sections}
\crefname{algocf}{Algorithm}{Algorithms}
\crefname{appendix}{Appendix}{Appendices}
\crefname{table}{Table}{Tables}
\crefname{problem}{Problem}{Problems}
\crefname{remark}{Remark}{Remarks}
\crefname{subsection}{Section}{Sections}
\crefname{subsubsection}{Section}{Sections}

\crefformat{equation}{(#2#1#3)}
\Crefformat{equation}{(#2#1#3)}

\renewenvironment{problem}
{\begin{mdframed}[topline=true,
                  bottomline=true,
                  leftline=false,
                  rightline=false,
                  linewidth=0.5pt]%
  \begin{probleminner}}
{\end{probleminner}\end{mdframed}}

\usepackage{accents}

\definecolor{aog}{rgb}{0.0, 0.5, 0.0}

\newcommand{\R}{\mathbb{R}}

\newcommand{\Z}{\mathbb{Z}}

\hypersetup{
  pdftitle={A Geometrical Acoustics based Focusing Algorithm for Layered Media in Medical Ultrasound},
  pdfauthor={S. Hackl, S. Hubmer, and R. Ramlau}
}

\title{A Geometrical Acoustics based Focusing Algorithm for Layered Media in Medical Ultrasound\thanks{Author-accepted manuscript. Accepted for publication in \emph{SIAM Journal on Imaging Sciences}. \copyright~by SIAM.}}
\author{
Simon Hackl\footnote{Johannes Kepler University Linz, Institute of Industrial Mathematics, Altenbergerstra{\ss}e 69, A-4040 Linz, Austria, (simon.hackl@jku.at), Corresponding author.} ,
Simon Hubmer\footnote{Johannes Kepler University Linz, Institute of Industrial Mathematics, Altenbergerstra{\ss}e 69, A-4040 Linz, Austria, (simon.hubmer@jku.at)} ,
Ronny Ramlau\footnote{Johann Radon Institute for Computational and Applied Mathematics, Altenbergerstra{\ss}e 69, A-4040 Linz, Austria, (ronny.ramlau@ricam.oeaw.ac.at)}\,\,\textsuperscript{,}\footnote{Johannes Kepler University Linz, Institute of Industrial Mathematics, Altenbergerstra{\ss}e 69, A-4040 Linz, Austria, (ronny.ramlau@jku.at)}
}

\begin{document}

\maketitle

\begin{abstract}
Ultrasound imaging is a widely used, non-invasive diagnostic tool in modern medicine. A crucial assumption is a constant sound speed in the observed medium. For large scale sound speed variations, this assumption leads to blurred and distorted images. In this paper, we present a Geometrical Acoustics Focusing Algorithm (GEOFA) which is able to correct for these aberrations, given a known layered medium setting with continuously differentiable medium boundaries. Existence and uniqueness conditions for a solution to the underlying system of equations are given. Using numerical simulations, the precision of our method is evaluated. Finally, the resulting image quality improvements are demonstrated in a phantom-based experimental setup.

\smallskip
\noindent \textbf{Keywords.} Focused Ultrasound, Aberration Correction, Geometrical Acoustics, Snell's law, Time of Flight, Focusing Delays, Layered Media, Ray Tracing.

\smallskip
\noindent \textbf{MSC codes.} 65Z05, 68U10, 78A10
\end{abstract}


\section{Introduction}
\begin{figure}[ht!]
    \centering
    \includegraphics[width=.85\textwidth]{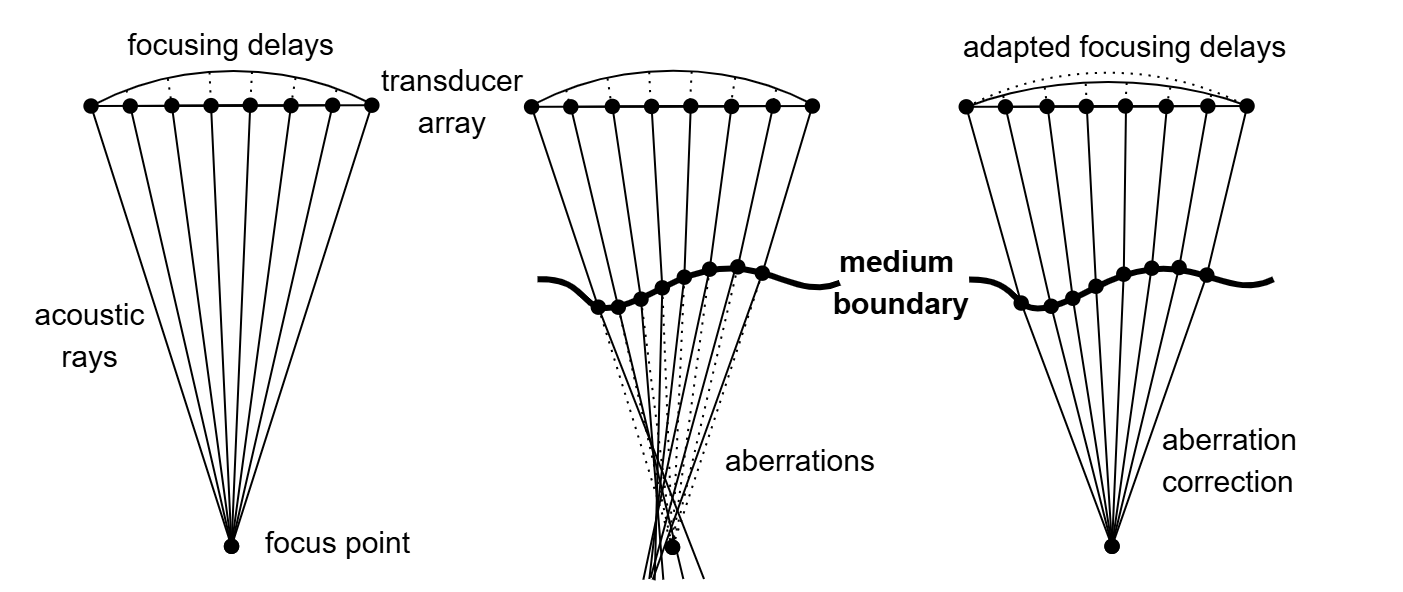}
    \caption{
    Left: Classic homogeneous medium focusing algorithm, Center: Aberrations in a multilayered medium, Right: Layer adapted aberration correction.}
    \label{fig:Distortions_of_homogeneous_foc_alg_in_multilayered_medium}
\end{figure}

Medical ultrasound imaging leverages the acoustic properties of biological tissues to provide real-time, non-invasive diagnostic capabilities. By precisely focusing ultrasound waves at specific depths within the patient's body, internal structures can be visualized with remarkable detail \cite{Overview_Ultrasound_Imaging}. This visualization is possible because sound waves are partially backscattered at tissue inhomogeneities along their propagation path. Backscattered waves, arriving at the ultrasound transducer, are recorded, and the ultrasound image is created in a procedure called \emph{Receive Focusing}. Due to the large amount of data involved in Receive Focusing, the simplifying assumption of a constant sound speed inside the observed medium is made, cf. \cref{fig:Distortions_of_homogeneous_foc_alg_in_multilayered_medium} and \cite{So_you_think_you_can_DAS,Overview_Ultrasound_Imaging}. The validity of this constant sound speed assumption depends on the setting \cite{Ultrasound_propagation_velocity_important_for_image_Quality,Imaging_quality_deteriorates_for_obese_patients, Acoustic_Impedance_Matching,TD_Cover_reconstruction_and_influence_on_image_quality}. The presence of large areas with a different sound speed creates various image distortions (\cite{ToF_Distortions_caused_by_Large_Scale_Tissue_layers, TD_Cover_reconstruction_and_influence_on_image_quality} and \cref{fig:Distortions_of_homogeneous_foc_alg_in_multilayered_medium}), collectively referred to as aberrations. Body parts containing fatty tissue or bone provide clinical examples of this issue \cite{Overview_of_Aberration_Correction_Methods}. Additionally, even the plastic cap on the ultrasound transducer may cause significant aberrations, cf. \cite{TD_Cover_reconstruction_and_influence_on_image_quality} and \cref{sect_4:numerical_simulations}.

In this paper, we introduce a geometrical acoustics-based aberration correction method which adapts focusing delays to a known layered medium via Time of Flight (ToF) calculations, as shown in \cref{fig:Distortions_of_homogeneous_foc_alg_in_multilayered_medium}. Historically, there have been various attempts for aberration correction \cite{Overview_of_Aberration_Correction_Methods}. In the area of enhanced precision in the ToF calculation, either the wave equation is solved \cite{liu1997propagation,mast1999simulation,mast1997simulation,pinton2011sources, shen2010computational} or straight ray tracing (SRT) methods \cite{Sound_Speed_Estimation_using_SAFA_model,mast1997simulation,vraalstad2024coherence} are used. A solution of the wave equation offers accuracy but high computational costs. In comparison, SRT is significantly faster, but inaccurate in cases of high sound speed variations, since refraction is neglected, see \cite{Overview_of_Aberration_Correction_Methods} and \cref{sect_4:numerical_simulations} of this work. 

An improvement to SRT can be achieved by additionally considering refraction, i.e., by applying geometrical acoustics. The approach most closely related to the one proposed in this paper is a geometrical acoustics based ToF calculation method for known layered media \cite{MINEO2022106747,MINEO2021106330}, which provides a satisfactory solution in its original context of material testing. However, a mathematical analysis of the underlying system of equations, including existence and uniqueness of solutions and the stability of its numerical solution, is not provided, and the approach has not been analyzed or validated in a medical ultrasound imaging setting. Moreover, since any ray tracing method requires the solution of a nonlinear system of equations for each transducer element-focus point pair, no such method, including ours, is currently suitable for real-time medical ultrasound imaging \cite{ibrahim2017efficient}.

The contributions of this paper are as follows. First, we establish theoretical and empirical foundations for focusing in known layered media, including existence and uniqueness conditions for the solution of the underlying problem. Second, we derive the Geometrical Acoustics Focusing Algorithm (GEOFA) and validate it both in k-Wave simulations and in a phantom-based experimental setting. The simulations show that GEOFA calculates ToFs and focusing delays with high accuracy, with errors typically below $1/32$ of the fundamental period. In the experiment, this precision translates to significant image quality improvements compared to the standard homogeneous medium focusing algorithm and to SRT. Our third contribution lies in demonstrating that GEOFA, in contrast to the method of \cite{MINEO2022106747,MINEO2021106330}, remains numerically stable up to the physical limit of total reflection. These contributions provide a theoretical and empirical basis for subsequent geometrical acoustics based algorithms with potential for real-time implementation in clinical ultrasound systems.

The outline of this paper is as follows: In \cref{sect_background}, we review basic medical ultrasound terminology, discuss the currently used homogeneous medium focusing algorithm and the SRT method, and explain geometrical acoustics. \cref{sect:3_mathematical_content} derives the system of equations our method is based on. Additionally, the problem is analyzed mathematically by proving existence and uniqueness conditions for the solution of the underlying system of equations. In \cref{sect_4:numerical_simulations}, our method is evaluated through numerical simulations and compared to a similar algorithm from \cite{MINEO2022106747}. Finally, a phantom-based experiment demonstrates our method's ability to remove the aberrations caused by a layered medium.

\section{Medical Ultrasound and Geometrical Acoustics}\label{sect_background}

\subsection{The Basic Setting}

In medical ultrasound imaging, sources of individual acoustic pressure waves, called transducer elements, are arranged in a transducer array inside the ultrasound transducer. The same transducer array receives backscattered acoustic signals from within the observed medium, which are used for image reconstruction \cite{Overview_Ultrasound_Imaging, szabo2013diagnostic}. Due to the high sound speed in medical settings ($\approx 1540 \, \text{m}/\text{s}$), multiple transmissions can be performed in a very short time period. Therefore, a static medium which does not change during exposure time is assumed. Finally, the standard coordinate system used throughout the ultrasound imaging literature is defined as follows \cite{szabo2013diagnostic}:

\begin{definition}[Ultrasound Imaging Coordinate system \cite{szabo2013diagnostic}]
    In the standard right-handed 3D coordinate system used in ultrasound imaging, the $x$-axis (\emph{azimuth} direction) and $z$-axis (\emph{depth} direction) span the primary imaging plane, with the $x$-axis corresponding to the horizontal and $z$-axis to the vertical image direction. The $y$-axis (\emph{elevation} direction) is perpendicular to this imaging plane.
\end{definition}

\subsection{Standard Focusing Method used in Homogeneous Media} \label{sect:2_2_one_layered_focusing_algorithm}

\begin{figure}[ht!]
    \begin{minipage}{.475\textwidth}
    
        \centering
        \includegraphics[width=\textwidth]{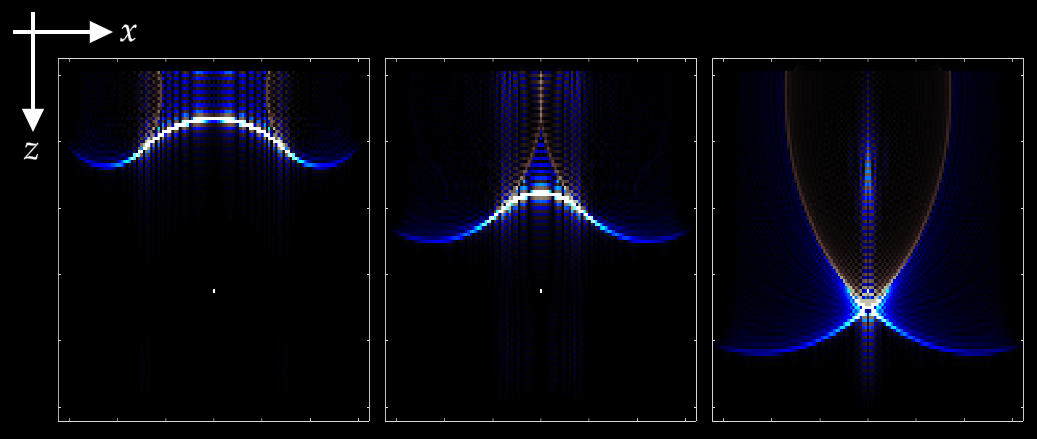}
    
    \end{minipage}
    \hfill
    \begin{minipage}{.49\textwidth}
    
        \centering
        \includegraphics[width=\textwidth]{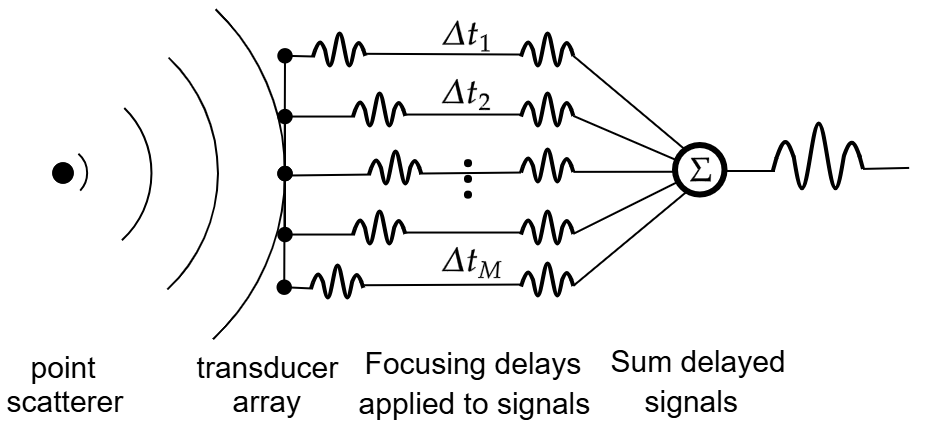}

    \end{minipage}

    \caption{Left: Transmit Focused wave at multiple points in time, simulated using k-Wave. Right: The DAS-Algorithm: A single backscattered wave arrives at the transducer elements at different times, which is corrected using focusing delays before the summation of the transducer signals.}
    \label{fig:Receive_Focusing_via_DAS}
\end{figure}

In medical ultrasound imaging, the most widely used imaging technique is Delay-and-Sum (DAS) beamforming, which uses focusing in two important steps \cite{So_you_think_you_can_DAS, Overview_Ultrasound_Imaging}. First, a \textit{Transmit Focus} is formed by delaying the emission of the individual transducer elements in time, as shown on the left of \cref{fig:Receive_Focusing_via_DAS} and in \cref{fig:Distortions_of_homogeneous_foc_alg_in_multilayered_medium}. These time delays compensate for the varying ToFs from the transducer elements to the desired transmit focus point, ensuring a simultaneous arrival of all emitted waves to achieve constructive interference. Inside the body, the transmitted sound wave is partially scattered at structures with a varying acoustic impedance $Z = c \rho$, where $c$ is the sound speed and $\rho$ the density of the tissue \cite{Wave_Prop_Time_Rev_Randomly_Layered_Media}. Reflections particularly happen at medium interfaces like skin-fat, fat-muscle, soft tissue-bone, or feto-maternal tissue boundaries. Backscattered acoustic signals detected by the transducer array are combined using \textit{Receive Focusing} to reconstruct the final ultrasound image. To be more precise, during Receive Focusing, appropriate time delays are applied to signals from each element before summation, representing the second critical focusing step in delay-and-sum beamforming, cf. the right of \cref{fig:Receive_Focusing_via_DAS}. This process is repeated for each image pixel, causing a significant computational effort \cite{Thomenius_Evolution_Beamforming}. 

A definition of focusing delays first requires a precise ToF definition.
\begin{definition}[ToF] \label{def:Time_of_Flight}
    The ToF between two points, called the source and focus points, is defined as the minimum time required for a spherical acoustic wave emitted at the source point to reach the focus point.
\end{definition}
\begin{remark}
    The ToF is given by Fermat's principle of least time. In layered media, it can be calculated explicitly using geometrical acoustics, see \cref{sect:2_3_geometrical_acoustics}.
\end{remark}

Time delays, applied to the transducer elements during transmit and receive focusing, are called \textit{focusing delays}, cf. \cite{Thomenius_Evolution_Beamforming}, \cref{fig:Distortions_of_homogeneous_foc_alg_in_multilayered_medium} and Definitions \ref{def:Transmit_Focusing} and \ref{def:acoustic_lens}. Transmit Focusing Delays compensate for varying ToFs from individual transducer elements to the desired transmit focus point, ensuring simultaneous arrival and constructive interference of all acoustic waves at the transmit focus point.
\begin{definition}[Transmit Focusing] \label{def:Transmit_Focusing}
    Let $(t_m)_{m = 1}^M$ be ToFs from acoustic\\ sources $(P_{0,m})_{m = 1}^M$ to a transmit focus point $P_N$. The associated \emph{transmit focusing delays} $(\Delta t_{m}^{Tx})_{m = 1}^M$ are defined as
    \begin{equation} \label{eq:acoustic_lens}
        \Delta t_{m}^{Tx} := \max \{ t_j | j = 1,\ldots,M \} - t_m \, .
    \end{equation}
    Furthermore, let $P_R$ be another point, $t_{R,N}$ the ToF from $P_R$ to $P_N$, and the sign $s$ be $1$ if the transmit focused wave reaches $P_N$ before $P_R$ and $-1$ otherwise. Then, the \emph{ToF of a wave transmit focused on $P_N$} to a point $P_R$ is defined as
    \begin{equation} \label{def:Transmit_ToF_to_RF_point}
        t_{\operatorname{Transmit}}(P_R) = \max \{ t_j | j = 1,\ldots,M \} + s \cdot t_{R,N} \, .
    \end{equation}
\end{definition}

Receive Focusing Delays account for two distinct propagation paths. First, the ToFs from the receive focus point back to the transducer elements have to be considered, cf. \cref{fig:Receive_Focusing_via_DAS} (right). Second, calculating the ToF that the transmitted wave takes to reach the receive focus point requires a model for the wave-front propagation during transmit, see \cref{def:Transmit_ToF_to_RF_point} and \cite{virtualSourceSAFT}.

\begin{definition}[Receive Focusing Delays] \label{def:acoustic_lens}
Let $(t_m)_{m = 1}^M$ be ToFs from acoustic sources $(P_{0,m})_{m = 1}^M$ to a receive focus point $P_R$. Furthermore, let $t_{\operatorname{Transmit}}(P_R)$ be the ToF of a transmit focused wave to $P_R$. Then, the \emph{receive focusing delays} $(\Delta t_{m}^{Rx})_{m = 1}^M$ are defined as
    \begin{equation} \label{eq:acoustic_lens_receive}
        \Delta t_{m}^{Rx} := t_{\operatorname{Transmit}}(P_R) + t_m \, .
    \end{equation}
\end{definition}

This paper presents an improved method for calculating the ToFs $(t_m)_{m = 1}^M$, which can subsequently be used to determine focusing delays, see (\ref{eq:acoustic_lens}), (\ref{def:Transmit_ToF_to_RF_point}) and (\ref{eq:acoustic_lens_receive}). Both transmit and receive focusing require ToF calculations from elements to focus points. For simplicity, we therefore do not distinguish between transmit and receive focusing in this paper.

The currently used Homogeneous Medium Focusing Algorithm (HMFA, \cite{So_you_think_you_can_DAS}) consists of two steps: First, assuming a homogeneous medium with a constant, predefined sound speed $c$, the ToFs $(t_m)_{m = 1}^M$ from sources $(P_{0,m})_{m = 1}^M$ to a focus point $P_1$ are calculated. Finally, the focusing delays are calculated, see Definition \ref{def:acoustic_lens}. This is summarized in \cref{Algo:homogeneous_medium_focusing}.

\begin{algorithm}[ht!]
\caption{Homogeneous Medium Focusing Algorithm (HMFA)}
\label{Algo:homogeneous_medium_focusing}
\LinesNumbered

\For{ \text{each source position} $P_{0,m}$ ($m = 1, \ldots, M$)}{
    \texttt{calculate ToF } $t_m = | P_1 - P_{0,m} | / c$
}
\For{ \text{each source position} $P_{0,m}$ ($m = 1, \ldots, M$)}{
        \texttt{calculate focusing delays} (\ref{eq:acoustic_lens}), (\ref{eq:acoustic_lens_receive})
}
\end{algorithm}

HMFA is implemented in ultrasound imaging devices with an assumed constant sound speed of $c = 1540 \, \text{m}/\text{s}$ \cite{Ultrasound_propagation_velocity_important_for_image_Quality,Thomenius_Evolution_Beamforming}. Methods for adapting this constant sound speed to the underlying medium using a global sound speed estimate have been suggested \cite{SoS_estimate_steering}, which are suitable for media with moderate sound speed contrast. However, for media containing areas with significantly deviating sound speeds, any calculations assuming a constant sound speed lead to incorrect ToFs at certain depths, which yields inaccurate focusing delays and consequently degraded image quality \cite{ToF_Distortions_caused_by_Large_Scale_Tissue_layers, TD_Cover_reconstruction_and_influence_on_image_quality}. We demonstrate this effect in \cref{sect_4_3_fat_layer_image_improvement} by comparing HMFA images obtained with a range of different constant sound speeds to the method suggested in this paper.

Focusing delay accuracy is of equal importance for non-imaging applications such as high-intensity focused ultrasound (HIFU) \cite{HIFU_Focus_point_Displacement} or material testing \cite{MINEO2022106747,MINEO2021106330}. By calculating the correct ToFs in layered media, we aim to correct these errors and thus expect improved imaging results. Clinically, a notable instance of sound speed deviations with significant impact on image quality occurs in ultrasound imaging of obese patients, where tissue heterogeneity causes substantial image degradation \cite{Ultrasound_propagation_velocity_important_for_image_Quality,Imaging_quality_deteriorates_for_obese_patients}. Moreover, already the plastic cap on top of the ultrasound transducer array, added to protect the transducer elements and to efficiently couple the transducer to the underlying tissue, causes significant aberrations, cf. \cite{Acoustic_Impedance_Matching,TD_Cover_reconstruction_and_influence_on_image_quality} and \cref{sect_4:numerical_simulations}. In this paper, we introduce an  algorithm which corrects these aberrations by calculating adapted focusing delays, as outlined on the right panel of \cref{fig:Distortions_of_homogeneous_foc_alg_in_multilayered_medium}.

\subsection{Geometrical Acoustics} \label{sect:2_3_geometrical_acoustics}

In order to improve the ToF calculation, we utilize the model of geometrical acoustics, which describes sound propagation using acoustic rays that are orthogonal to the wavefront. Geometrical acoustics is a frequently used approach \cite{Sound_Speed_Estimation_using_SAFA_model, Refraction_Based_SoS_estimation, Jaeger_2022, Geometrical_Acoustics_Derivation, Geometrical_Acoustics_for_solution_of_visualization_problems, MINEO2022106747, MINEO2021106330, Odegaard1995, Multilayer_Synthetic_Aperature_Focusing, TD_Cover_reconstruction_and_influence_on_image_quality}.

\begin{figure}[ht!]
    \centering
    \includegraphics[width=.95\textwidth]{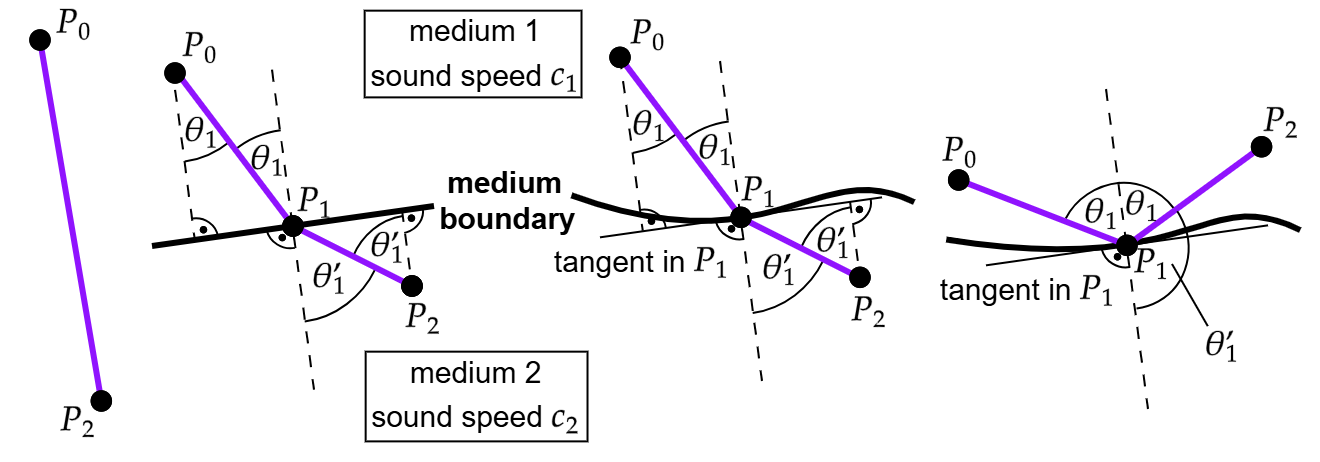}
    \caption{Geometrical acoustics is used in the following settings: Homogeneous medium (Left), Layered medium with straight boundary (Center left). Layered medium with $C^1$ medium boundary in case of transmission (Center right) or reflection (Right).}
    \label{fig:geometrical_acoustics}
\end{figure}

In this paper, geometrical acoustics is used to calculate the ToF through layered media, see \cref{fig:geometrical_acoustics}. In a homogeneous medium, the geometrical acoustics path, also called acoustic ray, is a straight line. On medium boundaries, refraction occurs, which means that the acoustic ray refracts according to Snell's law \cite{FeynmanLectureSnell}:
    \begin{equation} \label{eq:Snells_Law}
            c_2 \sin \theta_1 = c_1 \sin \theta_1' \, ,
    \end{equation}
where the angles $\theta_1$ and $\theta_1'$ are shown on the center left and center right panels of \cref{fig:geometrical_acoustics}. Another phenomenon occurring on medium boundaries is reflection. On any interface where the speed of sound changes, part of the acoustic wave is reflected. For reflections, the pre- and post-reflection angles are the same, as shown on the right panel of \cref{fig:geometrical_acoustics}. In cases of total reflection, the entire wave is reflected and no sound penetrates the medium boundary. Total reflection can be derived mathematically as follows: By rearranging Snell's law, we obtain $\sin \theta_1' = \frac{c_2}{c_1} \sin \theta_1$. If $\frac{c_2}{c_1} \sin \theta_1 \notin [-1,1]$, this equation cannot be fulfilled for any $\theta_1'$. In this case, total reflection occurs, i.e., $\theta_1' = \pi - \theta_1$, as outlined on the right panel of \cref{fig:geometrical_acoustics}.

Although geometrical acoustics is useful for ToF calculations, this approximate method has its limitations. E.g., it cannot model wave phenomena that arise from finite wavelength effects, such as diffraction, interference, or beam divergence \cite{geometrical_room_acoustic_modeling}.


\section{An Adapted Geometrical Acoustics Model for Focusing in Multilayered Media}\label{sect:3_mathematical_content}

\begin{figure}
    \centering
    \includegraphics[width=\textwidth]{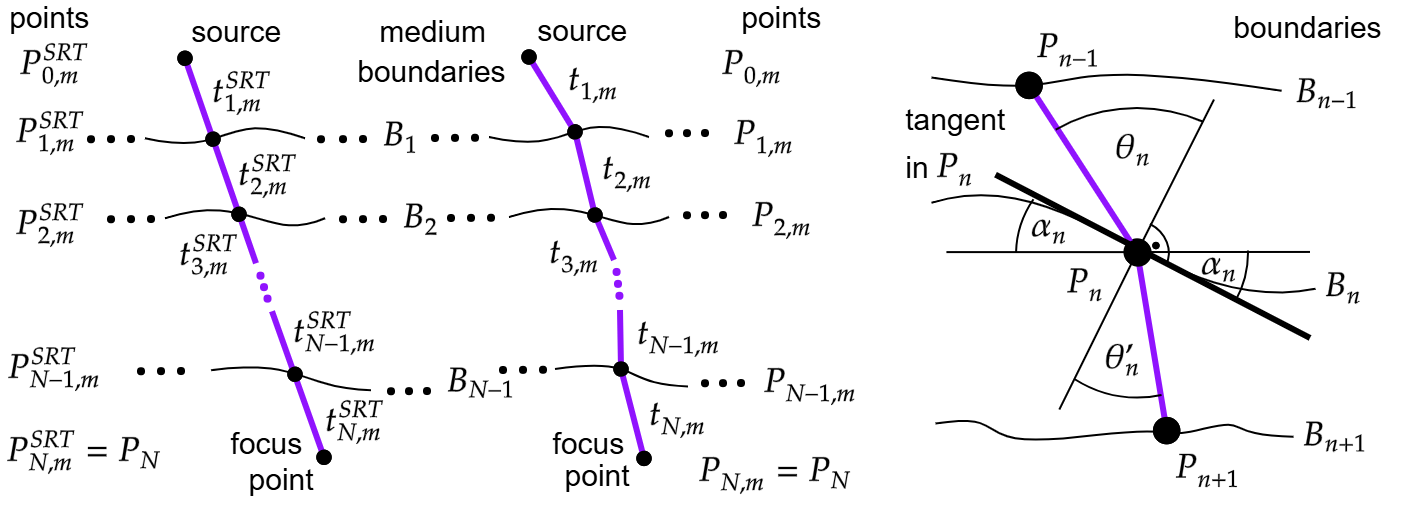}
    \caption{Both for SRT (Left) and GEOFA (Center), the ToF from the source $P_{0,m}$ to the focus point $P_N$ is calculated as the sum of the individual ToFs $(t_{n,m}^{\operatorname{SRT}})$ or $(t_{n,m})$, respectively, inside the layered medium. Right: In GEOFA, the medium boundary points are chosen s.t. the angles $\theta_n$ and $\theta_n'$ fulfill Snell's law of refraction (\ref{eq:Snells_Law}) in $P_n$.}
    \label{fig:SAFA_GOAT_Refraction}
\end{figure}

In this section we introduce our Geometrical Acoustics Focusing Algorithm (GEOFA) for multilayered media. \cref{sect:3_2_Defining_problem_and_notation} introduces the notation used. Furthermore, the problem of focusing waves in multilayered media is defined. Subsequently, in \cref{sect:3_4_MAFA}, we discuss GEOFA, which solves the multilayered focusing problem and improves HMFA by additionally considering the different sound speeds in the medium layers and refraction on the medium boundaries. In \cref{sect:3_3_MAFA_Existence_and_Uniqueness}, we derive conditions on the layered medium ensuring existence and uniqueness of solutions to the algorithm's underlying system of equations. Throughout the rest of this paper, we use the following assumptions:

\begin{assumption}[Known Medium Composition]\label{Ass:known_medium_composition_assumption}
    The medium is given by $N \geq 2$ layers, each having a constant sound speed $(c_n)_{n = 1}^N$. Both the sound speeds $(c_n)_{n = 1}^N$ and the medium boundaries $(B_n)_{n = 1}^{N-1}$ are known, the latter being of the form
    \[
        (B_n)_{n = 1}^{N-1} = ( \{ (x,b_n(x)) \, | \, x \in [0,1]\} )_{n = 1}^{N-1}\, ,
    \]
    where $b_n \in C^1([0,1])$ determines the medium boundary shape, see \cref{fig:SAFA_GOAT_Refraction}.
\end{assumption}

\begin{remark}
    The above assumption is justified by previous findings, e.g., in \cite{ToF_Distortions_caused_by_Large_Scale_Tissue_layers}, ``... most major arrival time fluctuations are caused by propagation through large-scale inhomogeneities such as fatty regions within muscle layers ...''.
\end{remark}    

\begin{remark}
    The assumption that the $z$-coordinate is a function of the $x$-coor\-dinate usually is fulfilled in applications. If this does not hold, one coordinate of the medium boundary can, at least locally, always be expressed as a function of the other coordinate. This assumption considerably simplifies the description of our algorithms compared to the more general setting of an arbitrary $C^1$ medium boundary.
\end{remark} 

Assumption \ref{Ass:known_medium_composition_assumption} is essential for the focusing algorithm presented in this paper, since our geometrical acoustics based ToF calculation approach relies on the underlying layered medium with $C^1$ medium boundaries being known. An example of a medium fulfilling this assumption is shown on the left of \cref{fig:SAFA_GOAT_Refraction}. 

Furthermore, we assume that the ToF in a layered medium can be calculated using geometrical acoustics, which is discussed in \cref{sect:2_3_geometrical_acoustics}:
    \begin{assumption}[Geometrical Acoustics] \label{Ass:Geometrical_acoustics}
        In layered media, Geometrical Acoustics models acoustic wave propagation using rays that are refracted at medium boundaries. Therefore, Geometrical Acoustics enables ToF calculation between two points.
    \end{assumption}
    \begin{remark}
        In addition to  multiple references in the literature confirming this assumption, such as \cite{ Sound_Speed_Estimation_using_SAFA_model, Refraction_Based_SoS_estimation, Jaeger_2022, Geometrical_Acoustics_Derivation, Geometrical_Acoustics_for_solution_of_visualization_problems, MINEO2022106747, MINEO2021106330, Multilayer_Synthetic_Aperature_Focusing}, we test it numerically in \cref{sect_4_1_ToF_measurements}, where we use ToF simulations to show that our model corrects the ToF errors of HMFA (\cref{Algo:homogeneous_medium_focusing}) occurring in multilayered media almost entirely. 
    \end{remark}

\begin{assumption}[Existence of Unique Medium Boundary Intersections]\label{Ass:unique_intersection}
    Each acoustic ray intersects each medium boundary $B_n$ exactly once.
\end{assumption}
\begin{remark}
    The above assumption simplifies the derivations in this section by assigning exactly one point per medium boundary to each acoustic ray. Note that a more general method allowing multiple intersections per boundary is possible using the same analytical framework considered here. In \cref{sect:3_3_MAFA_Existence_and_Uniqueness}, we derive sufficient conditions on the medium boundaries ensuring that this assumption holds.
\end{remark}

\subsection{Definitions and Problem Formulation} \label{sect:3_2_Defining_problem_and_notation}

In this section, we define the mathematical notation used throughout this paper and pose the multilayered focusing problem, which we solve below.

\begin{definition} \label{def:definition_of_notation}
    Let $M \geq 2$ acoustic sources be located in the positions ~\\$(P_{0,m})_{m = 1}^M = ( (x_{0,m}, z_{0,m}) )_{m = 1}^M$. 
    In the known layered medium setting (Assumption \ref{Ass:known_medium_composition_assumption}), define the \emph{medium boundary points} $(P_{n,m})_{n = 1,\, m = 1}^{N-1,\,M}$ located on the medium boundaries
    $( B_n )_{n = 1 }^{N-1}$ as
    \[
        P_{n,m} = (x_{n,m}, z_{n,m}) = (x_{n,m}, b_n(x_{n,m})) \, .
    \]
    Additionally define $(x_{N,m},z_{N,m}) = P_{N,m} := P_N$ for all $m = 1,\ldots,M$.    
\end{definition}

\begin{remark}
    In medical ultrasound imaging, the source positions $(P_{0,m})_{m = 1}^M$ and the focus point $P_N$ are given. In contrast, the medium boundary points $(P_{n,m})$ in between are unknown and have to be computed.
\end{remark}
The setting being fixed, the problem we aim to solve in this paper can be defined as follows.

\begin{problem}[Multilayered Focusing Problem] \label{prob:multilayered_focusing}
\begin{tabbing}
    \\ \textbf{Input}: \quad \= Medium properties: \quad \= Number of layers $N$,\\
                   \>                    \> Sound speeds $(c_n)_{n = 1}^N$ and\\
                   \>                    \> Medium boundaries $( B_n )_{n = 1}^{N-1}$.\\
                   \> Sources:           \> Number of sources $M$,\\
                   \>                    \> Position of the sources $(P_{0,m})_{m = 1}^M$.\\ 
                   \> Destination position: \> Focus point $P_N$.\\
    \\
    \textbf{Output}:        \> Focusing delays (\ref{eq:acoustic_lens}), (\ref{eq:acoustic_lens_receive}) for the sources $(P_{0,m})_{m = 1}^M$ to focus on $P_N$.
\end{tabbing}    
\end{problem}

\subsection{Time of Flight Calculation for given medium boundary points}\label{sect:3_2_ToF_calculation}

Assume that the medium boundary points are given, then the ToF can be calculated layer by layer as shown on the left of \cref{fig:SAFA_GOAT_Refraction}. To be more precise, let a fixed source $P_{0,m}$, medium boundary points $(P_{n,m})_{n = 1}^{N-1}$ and a focus point $P_{N}$ be given. The acoustic ray inside each layer is a straight line, because all medium layers are assumed to be homogeneous. Since the sound speeds $(c_n)_{n=1}^N$ are known, the ToFs $(t_{m,n})_{n=1}^N$ in each layer are given by $t_{n,m} = | P_{n,m} - P_{n-1,m}| /c_n$. The acoustic ray from $P_{0,m}$ to $P_N$ is a connection of all straight acoustic rays in each layer, cf. \cref{fig:SAFA_GOAT_Refraction}. Therefore, the total ToF $t_m$ from $P_{0,m}$ to $P_N$ is the sum of the ToFs in each layer, i.e.,
\begin{equation} \label{eq:ToF_calculation_for_given_medium_boundary_points}    
    t_m = \sum_{n = 1}^N t_{n,m} = \sum_{n = 1}^N \frac{ | P_{n,m} - P_{n-1,m} | }{c_n} \, .
\end{equation}
By performing this calculation for all sources $(P_{0,m})_{m = 1}^M$, all ToFs $(t_m)_{m = 1}^M$ are obtained, which can be used to calculate the focusing delays, cf. Definition \ref{def:acoustic_lens}. 

\subsection{Straight Ray Tracing (SRT)} \label{sect:3_3_SRT}

SRT is the current state-of-the-art method used for sound speed reconstruction and imaging in a known layered medium setting \cite{jaeger2015computed,sanabria2018spatial,vraalstad2024coherence}. For a fixed medium boundary $B_n$, source $P_{0,m}$ and focus point $P_N$, SRT determines the medium boundary points $P_{n,m}$ by intersecting the straight line through $P_{0,m}$ and $P_N$ with the medium boundary $B_n$, see the left part of \cref{fig:SAFA_GOAT_Refraction}. Using the medium boundary points, the ToF from $P_{0,m}$ to $P_N$ is subsequently calculated layer-wise, as discussed in \cref{sect:3_2_ToF_calculation} and focusing delays are derived from this ToF, cf.~Definition~\ref{def:acoustic_lens}.

\subsection{Geometrical Acoustics Focusing Algorithm (GEOFA)} \label{sect:3_4_MAFA}

In this section, a geometrical acoustics based focusing algorithm, which solves the Multilayered Focusing Problem \ref{prob:multilayered_focusing}, is formulated. Based on Geometrical Acoustics, we derive a system of equations for the unknown medium boundary points required to calculate the ToF through the medium, see \cref{sect:3_2_ToF_calculation}. For this derivation, we fix the index $m$ of the source and omit it for simplification, i.e., for fixed $m$ we write $(P_{n})_{n = 0}^{N-1}$ instead of $(P_{n,m})_{n = 0}^{N-1}$. Furthermore, we fix the considered medium boundary, denoted with the index $n$. Therefore, the following analysis is done for a point $(x_n,z_n) = P_n = P_{n,m}$ on a fixed  medium boundary $B_n$. 

In Snell's law, the pre- and post-refraction angles $\theta_{n}$ and $\theta'_{n}$ are crucial. Therefore, we include them into the set of unknowns in addition to $(x_n,z_n) = P_n$. The underlying geometry is shown on the right panel of \cref{fig:SAFA_GOAT_Refraction}. Additionally, the tangent slope of the medium boundary in $P_{n}$ has to be taken into account using the angle $\alpha_{n}$. Thus, for each medium boundary $B_{n}$, the five unknown variables, $x_{n}, z_{n}, \theta_{n}, \theta'_{n}$, and $\alpha_{n}$, have to be determined using the five equations we derive below.

First, refraction at the medium boundary is modeled using Snell's law (\ref{eq:Snells_Law}). Furthermore, the point $P_n$ is located at the medium boundary $B_n$, which yields
    \begin{equation} \label{eq:medium_boundary_equation}
        z_n = b_n(x_n) \, .
    \end{equation}
Moreover, the angle $\alpha_n$ of the tangent to the medium boundary in $P_n$ is given by
    \begin{equation}\label{eq:tangent_of_medium_boundary_equation}
        \tan \alpha_n = b_n'(x_n).
    \end{equation}
Up to this point, we have obtained three equations (\ref{eq:Snells_Law}), (\ref{eq:medium_boundary_equation}) and (\ref{eq:tangent_of_medium_boundary_equation}). In the following lemma, we derive two more equations for the terms $\sin \theta_n$ and $\sin \theta_n'$ which appear in Snell's law.
    
    \begin{figure}[ht!]
        \centering
        \includegraphics[width=.9\textwidth]{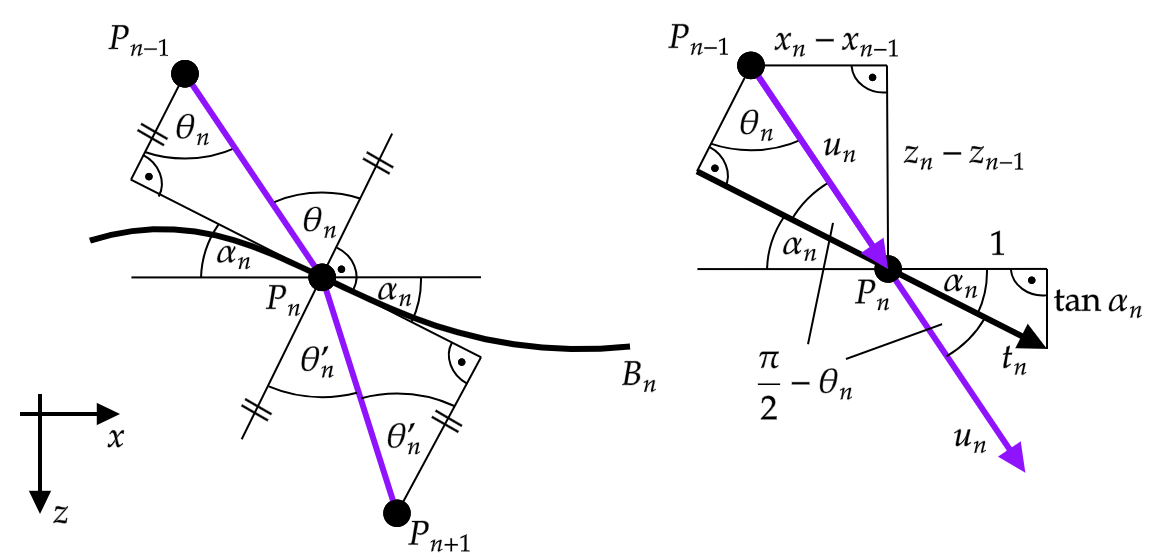}
        \caption{Left: Geometric construction in \cref{lemma:direct_equations_for_sin_theta}. Right: Defining vectors $u_n$ and $t_n$ and finding their enclosed angle.}
        \label{fig:mafa_geometric_analysis_of_refraction}
    \end{figure}

\begin{lemma} \label{lemma:direct_equations_for_sin_theta}
    For the coordinates of the points $P_{n-1}$, $P_n$ and $P_{n+1}$ and the angles $\alpha_n$, $\theta_n$ and $\theta'_n$, outlined in \cref{fig:mafa_geometric_analysis_of_refraction}, the following equations hold:
    \begin{align}
        \sin \theta_n &= \frac{x_{n} - x_{n-1} +( z_{n} - z_{n-1} )\tan \alpha_n}{ \sqrt{1 + \tan^2 \alpha_n} \sqrt{(x_{n} - x_{n-1})^2 + (z_{n} - z_{n-1})^2}} \, , \label{eq:first_medium_boundary_eq}\\
        \sin \theta'_n &= \frac{x_{n+1} - x_{n} +  (z_{n+1} - z_{n} ) \tan \alpha_n}{\sqrt{ 1 + \tan^2 \alpha_n } \sqrt{(x_{n+1} - x_{n})^2 + (z_{n+1} - z_{n})^2} } \, . \label{eq:second_medium_boundary_equation}
    \end{align}
\end{lemma}
\begin{proof}

    First, define $\Delta x_n := x_{n} - x_{n-1}$, $\Delta z_n := z_{n} - z_{n-1}$, and the vectors $t_n = (1, \tan \alpha_n)$ and $u_n = (\Delta x_n, \Delta z_n)$, which are illustrated in \cref{fig:mafa_geometric_analysis_of_refraction}. By adding the vector $u_n$ to $P_n$, the right panel of \cref{fig:mafa_geometric_analysis_of_refraction} is obtained. From this, it follows that
    \[
    \begin{aligned}
        \cos\bigg(\frac{\pi}{2} - \theta_n \bigg) &= \langle u_n, t_n \rangle \big/ ( \| u_n \|_2 \, \| t_n \|_2 ) 
        = \frac{\Delta x_n + \tan \alpha_n \Delta z_n }{\sqrt{1 + \tan^2 \alpha_n } \sqrt{ (\Delta x_n)^2 + (\Delta z_n)^2}} \, .
    \end{aligned}
    \]
    Due to $\sin(\theta_n) = \cos\big(\frac{\pi}{2} - \theta_n\big)$, (\ref{eq:first_medium_boundary_eq}) follows. Equation \ref{eq:second_medium_boundary_equation} can be derived analogously.
\end{proof}

\begin{definition}[GEOFA system of equations]
    The system of equations which consist of (\ref{eq:Snells_Law}), (\ref{eq:medium_boundary_equation}), (\ref{eq:tangent_of_medium_boundary_equation}), (\ref{eq:first_medium_boundary_eq}), and (\ref{eq:second_medium_boundary_equation}) for $n = 1,\ldots,N-1$, i.e.,
    \[
    \begin{aligned}
        c_{n+1} \sin \theta_n &= c_n \sin \theta_n' \, ,\\
        z_n &= b_n(x_n) \, ,\\
        \tan \alpha_n &= b_n'(x_n) \, , \\
        \sin \theta_n &= \frac{x_{n} - x_{n-1} +( z_{n} - z_{n-1} )\tan \alpha_n}{ \sqrt{1 + \tan^2 \alpha_n} \sqrt{(x_{n} - x_{n-1})^2 + (z_{n} - z_{n-1})^2}} \, , \\
        \sin \theta'_n &= \frac{x_{n+1} - x_{n} +  (z_{n+1} - z_{n} ) \tan \alpha_n}{\sqrt{ 1 + \tan^2 \alpha_n } \sqrt{(x_{n+1} - x_{n})^2 + (z_{n+1} - z_{n})^2} } \, ,
    \end{aligned}           
    \]
    is called \emph{GEOFA system of equations}.
\end{definition}

To summarize, we have derived the GEOFA system of equations, which consists of the $5\cdot(N-1)$ equations (\ref{eq:Snells_Law}), (\ref{eq:medium_boundary_equation}), (\ref{eq:tangent_of_medium_boundary_equation}), (\ref{eq:first_medium_boundary_eq}), and (\ref{eq:second_medium_boundary_equation}), for the $5\cdot(N-1)$ unknowns $x_n$, $z_n$, $\theta_n$, $\theta'_n$, and $\alpha_n$, where $n = 1,\ldots,N-1$. Therefore, to calculate the focusing delays required for an $N$ layered medium, a system of $5 \cdot (N-1)$ equations has to be solved to obtain the medium boundary points $(P_n)_{n = 1}^{N-1}$. From these, the ToF and the focusing delay can be calculated as outlined in Algorithm~\ref{Alg:Detailed_multilayered_focusing_Algorithm}.

\begin{algorithm}[ht!]
\caption{Geometrical Acoustics Focusing Algorithm (GEOFA)}
\label{Alg:Detailed_multilayered_focusing_Algorithm}
\LinesNumbered
\For{ each source position $P_{0,m}$ ($m = 1, \ldots, M$) }{
    \texttt{solve the GEOFA system of equations} $(\ref{eq:Snells_Law}),(\ref{eq:medium_boundary_equation}),(\ref{eq:tangent_of_medium_boundary_equation}), (\ref{eq:first_medium_boundary_eq}),(\ref{eq:second_medium_boundary_equation})$ \texttt{for the $5 \cdot (N-1)$ unknowns $x_n, z_n, \theta_n, \theta'_n, \alpha_n$, $n = 1, \ldots N-1$.}\\
    \texttt{calculate ToF:}
    \[
       t_m = \sum_{n=1}^{N} \sqrt{(x_{n}-x_{n-1})^2 + (z_n - z_{n-1})^2} \big/ c_n
    \]
}
\For{ \text{each source position} $P_{0,m}$ ($m = 1, \ldots, M$)}{
    \texttt{calculate focusing delays} (\ref{eq:acoustic_lens}), (\ref{eq:acoustic_lens_receive})
}
\end{algorithm}

For solving the nonlinear system of equations in Step 2 of \cref{Alg:Detailed_multilayered_focusing_Algorithm}, an iterative Newton-type method can be used. Due to moderate sound speed variations in medical ultrasound imaging, refraction is moderate. Hence, medium boundary points determined by SRT can be used as an initial guess for an iterative solution method.

\begin{remark}
    In \cref{sect_4_1_ToF_measurements}, we show that GEOFA corrects the ToF and focusing delay errors of HMFA and SRT below the level of significance for medical ultrasound imaging and thus solves Problem \ref{prob:multilayered_focusing}. Furthermore, in \cref{sect_4_3_fat_layer_image_improvement}, we demonstrate that in multilayered media, GEOFA corrects ultrasound image aberrations that occur in HMFA and SRT due to incorrect focusing delays.
\end{remark}

\paragraph{Comparison with a previously introduced method}

A different geometrical acoustics based method for multilayered media, which we refer to as MINEO, has been published in the context of non-destructive testing \cite{MINEO2022106747,MINEO2021106330}. Similar to GEOFA, MINEO determines the medium boundary points $(P_{n})_{n=1}^{N-1}$ by solving a nonlinear system of equations derived from Snell's law and the medium boundary geometry, and subsequently calculates the ToF using \eqref{eq:ToF_calculation_for_given_medium_boundary_points}. However, the two methods differ in their choice of unknowns, their formulation of the underlying equations, and consequently in the structure of the resulting system of equations.

MINEO uses the $x$-coordinates $(x_n)_{n=1}^{N-1}$ of the medium boundary points together with $N$ angles $(\gamma_n)_{n=1}^N$ as unknowns, which are derived by solving a system of $(2N-1)$ equations. Here, $\gamma_n$ denotes the angle between the acoustic ray in layer $n$ and the vertical axis, which, unlike the pre- and post-refraction angles $\theta_n$ and $\theta_n'$ in GEOFA, is not measured with respect to the boundary normal. Therefore, applying Snell's law at the boundary $B_n$ requires correcting $\gamma_n$ by the boundary tangent angle $\alpha_n = \arctan(b_n'(x_n))$, which yields the true pre-refraction angle $\theta_n = \gamma_n - \alpha_n$. In contrast, GEOFA explicitly uses the pre- and post-refraction angles $\theta_n$ and $\theta_n'$ as a separate variable, which can be inserted directly into Snell's law without any further correction.

As discussed in this section, GEOFA introduces, for each medium boundary $B_n$, the five unknowns $x_n$, $z_n$, $\theta_n$, $\theta_n'$ and $\alpha_n$, resulting in the formally larger GEOFA system of $5(N-1)$ equations. This system includes explicit structural equations such as $z_n = b_n(x_n)$ and $\tan\alpha_n = b_n'(x_n)$, which are not present in MINEO, where they are substituted into the remaining equations to obtain a smaller system. While this substitution reduces the number of unknowns, it increases the complexity of the individual equations, whereas GEOFA retains a larger but structurally simpler system.

In \cref{sect_4_2:ToF_errors_depending_on_med_bound_points}, we show that the GEOFA system of equations admits improved numerical stability compared to MINEO, with numerical solvers converging in all physical settings right until the physical limit of total reflection.

\subsection{Existence and Uniqueness Results} \label{sect:3_3_MAFA_Existence_and_Uniqueness}

In the first step of GEOFA (Line $2$ of \cref{Alg:Detailed_multilayered_focusing_Algorithm}), a nonlinear system of equations has to be solved. In this section, we therefore discuss conditions for the existence and uniqueness of solutions. Equivalently, the derived conditions ensure that Assumption \ref{Ass:unique_intersection} about the existence of unique medium boundary points is satisfied. Technical details can be found in \cref{sect:6_Appendix}.

\paragraph{Existence of a Solution in the GEOFA System of Equations} \label{sect:3_3_1_Existence_of_GOAT_solution}

\begin{figure}[ht!]
    \centering
    \includegraphics[width=\textwidth]{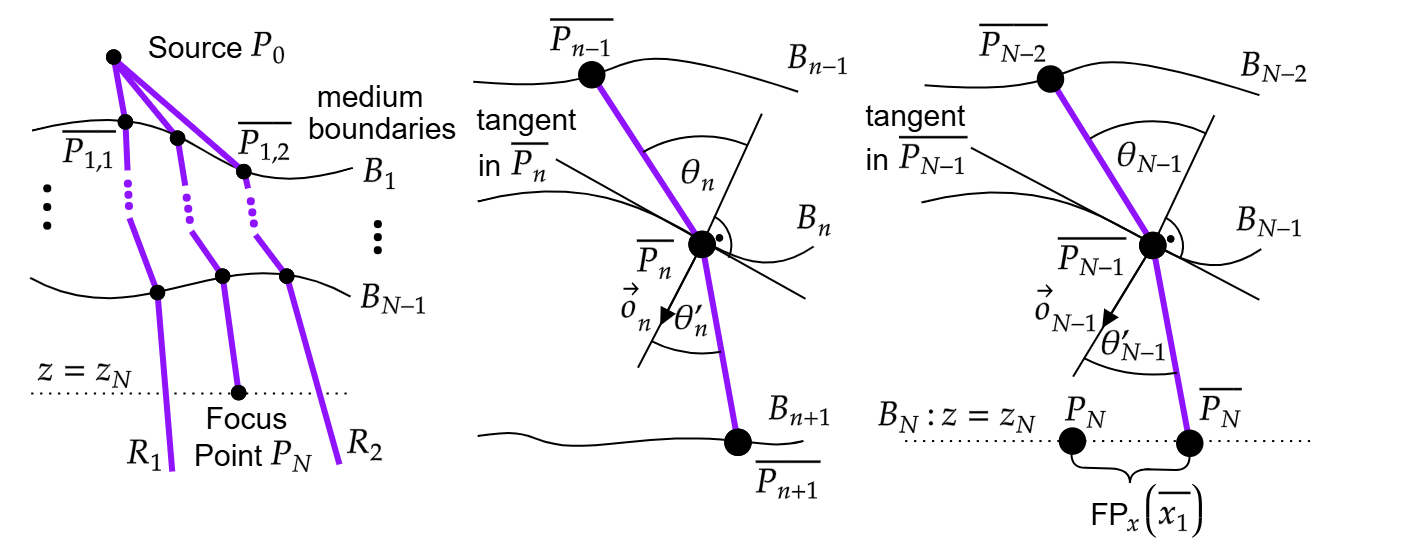}
    \caption{Left: Setting from \cref{thm:Existance_of_MAFA_solution}, Center: Refraction at the $n$-th medium boundary for $n = 1,\ldots,N-2$, Right: Refraction at the last medium boundary and definition of $\operatorname{FP}_x$.}
    \label{fig:exact_setting_in_MAFA_existence_proof}
\end{figure}

 \begin{definition}
    Let $R_1$, $R_2$ be two acoustic rays originating in $P_0$ and intersecting the medium boundary $B_1$ in the points $\overline{P_{1,1}}$ and $\overline{P_{1,2}}$, respectively. An acoustic ray originating in $P_0$ is defined to be \emph{between $R_1$ and $R_2$} if it intersects the first medium boundary in a point between $\overline{P_{1,1}}$ and $\overline{P_{1,2}}$, as shown on the left panel of \cref{fig:exact_setting_in_MAFA_existence_proof}.
\end{definition}

\begin{theorem}[Existence of a GEOFA solution]\label{thm:Existance_of_MAFA_solution}
    Assume that there exist two acoustic rays $R_1$ and $R_2$ originating in $P_0$, which pass the focus point $P_N$ on opposite sides, as indicated on the left panel of \cref{fig:exact_setting_in_MAFA_existence_proof}. Furthermore, assume that for all acoustic rays between $R_1$ and $R_2$, no total reflection and exactly one intersection with each medium boundary and with the line $B_N: z = z_N$ occurs. Then, there exists a solution for the GEOFA system of equations. 
\end{theorem}
    
\begin{proof}
    The underlying idea of this proof is to construct a continuous function $\operatorname{FP}_x$ whose roots are associated with a solution to the GEOFA system of equations. Then, the existence of a root of $\operatorname{FP}_x$ is shown using the intermediate value theorem.
    
    First, note that for an acoustic ray starting in $P_0$, the path through the medium is uniquely determined by its intersection with $B_1$ in $\overline{P_1}$. This is due to the underlying law of geometrical acoustics, which can be used to determine the acoustic ray's path through the medium once $P_0$ and $P_1$ are given. Next, define $\overline{P_N} = (\overline{x_N},\overline{z_N})$ as the intersection of the line $B_N: z = z_N$ with the acoustic ray originating in $P_0$ and passing through $\overline{P_1} = (\overline{x_1},b_1(\overline{x_1})) \in B_1$, as shown on the right of \cref{fig:exact_setting_in_MAFA_existence_proof}. Using $\overline{P_N}$, define $\operatorname{FP}_x(\overline{x_1}) = \overline{x_N} - x_N$. By determining the medium boundary points $(\overline{P_n})_{n=1}^{N-1}$ with the laws of geometrical acoustics, it can be shown that $\overline{x_N}$ continuously depends on $\overline{x_1}$, see \cref{lemma:medium_boundary_points_continuisly_dependence} in \cref{sect:6_1_Appendix_Existence_Result}. Therefore, $\operatorname{FP}_x(\overline{x_1})$ continuously depends on $\overline{x_1}$.

    Let $\overline{x_{1,1}}$ and $\overline{x_{1,2}}$ be the $x$-coordinates of the intersection points of $B_1$ with the acoustic rays $R_1$ and $R_2$. Due to $R_1$ and $R_2$ passing $P_N$ on opposite sides, $\operatorname{FP}_x(\overline{x_{1,1}})$ and $\operatorname{FP}_x(\overline{x_{1,2}})$ have different signs. Furthermore, the function $\operatorname{FP}_x$ is well defined on the interval $[\min\{\overline{x_{1,1}}, \overline{x_{1,2}}\}, \max\{ \overline{x_{1,1}},\overline{x_{1,2}} \}]$ due to the assumed properties of all acoustic rays between $R_1$ and $R_2$. Since $\operatorname{FP}_x$ is continuous, the intermediate value theorem yields the existence of $x_1$ such that $\operatorname{FP}_x(x_1) = 0$. By the definition of $\operatorname{FP}_x$, this means that the acoustic ray starting in $P_0$ and intersecting $B_1$ in $(x_1,b_1(x_1))$ hits the focus point $P_N$. This is equivalent to the existence of a solution to the GEOFA system of equations.
\end{proof}

\begin{remark}
    Certain medical ultrasound imaging settings can cause total reflection. An example for this would be tissue containing bone. HMFA does not consider total reflection, leading to significant image distortions. Therefore, the \textit{no total reflection} assumption is already used implicitly in HMFA. Similarly, the condition requiring the \textit{existence of acoustic rays passing the focus point on opposite sides} is only violated in settings involving total reflection. Finally, the \textit{unique intersection condition} is used as a simplifying assumption for the GEOFA algorithm. However, at higher computational complexity, it is possible to generalize GEOFA to additionally consider multiple intersections per medium boundary.
\end{remark}

\begin{remark}
    If the function $\operatorname{FP}_x$, defined in the proof of \cref{thm:Existance_of_MAFA_solution}, is strictly monotone, the uniqueness of the obtained solution additionally follows. Note that the monotonicity of $\operatorname{FP}_x$ depends on the given setting.
\end{remark}

\paragraph{Uniqueness Condition} \label{sect:3_3_2_Uniqueness_of_GOAT_solution}

In order to derive a uniqueness condition for the GEOFA system of equations, we analyze the ToF of all possible acoustic rays, including those that do not fulfill Snell's law. This analysis requires the following definitions.

\begin{notation}
    To abbreviate the difference in distance in $x$ and $z$ direction of two consecutive medium boundary points, we use the notation
    \[
    \begin{aligned}
       \Delta x_n &:= x_n - x_{n-1}, \qquad \qquad
       \Delta b_n(x_n) := b_n(x_n) - b_{n-1}(x_{n-1}) \, .
    \end{aligned}
    \]
    Here, we use the conventions $b_0(x_0) := z_0$ and $b_N(x_N) := z_N$.
\end{notation}
In contrast to \cref{def:Time_of_Flight}, where the ToF between two points was defined in the sense of Fermat's principle of least time, now we want to define the ToF along an arbitrary acoustic ray, which is a function of its medium boundary intersections. 
\begin{definition}
    The function ToF: $B_1 \times \ldots \times B_{N-1} \rightarrow \R^+$ denotes the \emph{Time of Flight} along the acoustic ray from $P_0$ to $P_N$ through $P_n = (x_n,b_n(x_n)) \in B_n,$ for all $n = 1,$ $\ldots,$ $N-1$:
    \[
        \text{ToF}(x_1,\ldots,x_{N-1}) := \sum_{n = 1}^N |P_n - P_{n-1}|/c_n = \sum_{n = 1}^N \sqrt{(\Delta x_n)^2 + (\Delta b_n(x_n))^2}/c_n \, .
    \] 
\end{definition}

Using these definitions, we can show that Snell's law being satisfied is equivalent to the ToF function having a stationary point, which is in turn equivalent to the slopes of the medium boundary satisfying a certain equation.

\begin{theorem} \label{thm:uniqueness_GOAT_TFAE_statement}
    For $1 \leq n \leq N-1$, let $(x_n,b_n(x_n)) = P_n \in B_n$. Then the following statements are equivalent as long as (\ref{eq:med_bound_slope_for_equal_ToF}) has a non-zero denominator. 
            \begin{enumerate}
                \item $(x_1,\ldots,x_{N-1})$ is a stationary point for the ToF-function.
                \item Snell's law (\ref{eq:Snells_Law}) holds in $P_1$, $\ldots$, $P_{N-1}$.
                \item For all $n = 1,\ldots,N-1$ the medium boundary slope $b_n'$ satisfies
    \begin{equation} \label{eq:med_bound_slope_for_equal_ToF}
        b_n'(x_n) = \frac{c_n | P_n - P_{n-1} | \Delta x_{n+1} - c_{n+1} | P_{n+1} - P_n | \Delta x_n }{c_{n+1} | P_{n+1} - P_n | \Delta b_n(x_n) - c_n | P_n - P_{n-1} | \Delta b_{n+1}(x_{n+1}) } \, . 
    \end{equation}
            \end{enumerate}
\end{theorem}
\begin{proof}
    See Appendix \ref{sect:6_2_Appendix_Uniqueness_Result}.
\end{proof}

\begin{remark}
    Snell's law being valid exactly for the stationary points of the ToF-function (equivalence of 1 and 2) is a special case of Fermat's principle, which is discussed in \cite{Chaves_nonimaging_optics}.
\end{remark}

\begin{corollary}[Uniqueness condition for the GEOFA system of equations] \label{thm:Uniqueness_GOAT}
    There exists a unique solution to the GEOFA system of equations if and only if on every medium boundary $B_n$, there exists a unique point $P_n =(x_n,b_n(x_n))$ fulfilling (\ref{eq:med_bound_slope_for_equal_ToF}).
\end{corollary}
\begin{proof}
    Assume that on every medium boundary $B_n$, there exists a unique point $P_n = (x_n,b_n(x_n))$ fulfilling (\ref{eq:med_bound_slope_for_equal_ToF}). In this case, \cref{thm:uniqueness_GOAT_TFAE_statement} yields that Snell's law holds in precisely these points. Thus, by inserting the coordinates of these points into the GEOFA system, the unknowns $x_n, z_n, \theta_n, \theta'_n, \alpha_n$ can be uniquely determined.

    Conversely, assume that there exists a unique solution to the GEOFA system of equations. This implies that at every medium boundary, there exists a unique point where Snell's law is fulfilled, because the other equations in the GEOFA system always hold if $P_n \in B_n$ for all $n = 1, \ldots, N-1$. Therefore, \cref{thm:uniqueness_GOAT_TFAE_statement} yields that at every medium boundary $B_n$, there is a unique point $P_n = (x_n,b_n(x_n))$ fulfilling (\ref{eq:med_bound_slope_for_equal_ToF}). 
\end{proof}

We conclude this section by examining the special case of $N = 2$ layers. There, \cref{thm:Uniqueness_GOAT} can be visualized. 

\begin{definition}
    For a given point $P'$, the set consisting of all points $P$ fulfilling $\text{ToF}(P) = \text{ToF}(P') $ is called the \emph{ToF level set}. Furthermore, define the \emph{equivalence relation} $P \sim P'$ for points $P, P'$ on the same ToF level set.
\end{definition}

\begin{figure}[ht!]
    \centering
    \includegraphics[width=0.425\linewidth]{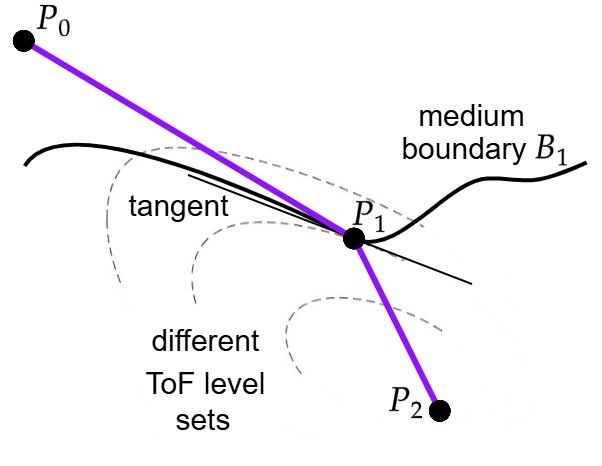}
    \caption{A 2-layered medium with a medium boundary $B_1$ and different ToF level sets (dashed). $P_1$ is the medium boundary point on the ToF level set with minimal ToF. In $P_1$, the tangents to the medium boundary and the ToF level set $[P_1]_\sim$ coincide.}
    \label{fig:Def_of_equal_ToF_lines_and_angle_to_real_boundary}
\end{figure}

\begin{remark}
    A Cartesian oval with respect to $Q_1, Q_2 \in \R^2$ and $a,b \in \R$ is defined as the set of all points $P$ s.t. $|P - Q_1| + a |P - Q_2| = b$ \cite{villarino2007decartesperfectlens}. Since the equation $\text{ToF}(P) = \text{const}$ can be rearranged into this form, ToF level sets are Cartesian ovals. Furthermore, by interpreting (\ref{eq:med_bound_slope_for_equal_ToF}) as an ODE with a given point $P$ as initial condition, the ToF level set $[P]_\sim$ can be computed numerically. Examples for ToF level sets in a 2-layered medium starting at different points are the dashed oval curves visualized in \cref{fig:Def_of_equal_ToF_lines_and_angle_to_real_boundary}.
\end{remark}

\begin{remark}
    The slope of the ToF level set depends on the previous and subsequent points, see (\ref{eq:med_bound_slope_for_equal_ToF}). Therefore, ToF level sets are most useful for 2-layered media, where $P_0$ and $P_2$, and thus the ToF level sets themselves are known a priori. For $N>2$ layers, the ToF level sets depend on a priori unknown neighboring medium boundary points.
\end{remark}

\begin{corollary}
    In a 2-layered medium, there exists a unique medium boundary point fulfilling Snell's law if and only if there exists a unique $P_1 \in B_1$ touching the ToF level set $[P_1]_\sim$.
\end{corollary}

\begin{proof}
    First, note that $P_1 \in B_1$ touching the ToF level set $[P_1]_\sim$ is equivalent to $P_1$ being a stationary point of the ToF-Function along $B_1$. Assuming there exists a unique medium boundary point $P_1$ fulfilling Snell's law, \cref{thm:uniqueness_GOAT_TFAE_statement} yields that $P_1$ is a stationary point of the ToF-Function, and thus touching the ToF level set $[P_1]_\sim$. Conversely, if there exists a  unique $P_1 \in B_1$ touching $[P_1]_\sim$, $P_1$ is the unique stationary point of the ToF-Function. Due to \cref{thm:uniqueness_GOAT_TFAE_statement}, we therefore obtain that $P_1$ is the unique medium boundary point fulfilling Snell's law.
\end{proof}

\subsection{Special Medium Boundary Shapes}

In this section, we simplify the GEOFA system of equations for the special cases of straight and elliptic medium boundaries. These settings are of particular interest, since they correspond to the cases of a fat layer above a tissue layer or an ultrasound transducer cover. 

\paragraph{Straight medium boundary}

For a straight medium boundary $B_n$, the medium boundary function is of the form $ b_n(x) = k_n x + d_n$. This simplifies the GEOFA system, with (\ref{eq:medium_boundary_equation}) and (\ref{eq:tangent_of_medium_boundary_equation}) now reading
\begin{align}
    z_n &= b_n(x_n) = k_n x_n + d_n \, , \label{straight_medium_boundary_equation}\\
    \tan(\alpha_n) &= b_n'(x_n) = k_n\, . \label{straight_medium_boundary_slope}   
\end{align}
Equation (\ref{straight_medium_boundary_slope}) implies that for straight medium boundaries, 
(\ref{eq:tangent_of_medium_boundary_equation}) is independent of $x_n$. Therefore, the variable $\alpha_n$ can be omitted, since $\alpha_n = \arctan(k_n)$ is known a priori. Additionally, (\ref{straight_medium_boundary_slope}) can be inserted into (\ref{eq:first_medium_boundary_eq}) and (\ref{eq:second_medium_boundary_equation}), reducing the number of equations from five to four.

Finally, we consider the most simple setting of a flat, axis-parallel medium boundary $B_n$. In this setting, the medium boundary function is constant, i.e., $b_n = d_n$. Due to (\ref{eq:medium_boundary_equation}), $z_n = d_n$ is known a priori. Furthermore, (\ref{eq:tangent_of_medium_boundary_equation}) reads $\tan \alpha_n = b_n'(x_n) = 0 $. Using these simplifications and inserting (\ref{eq:first_medium_boundary_eq}) and (\ref{eq:second_medium_boundary_equation}) into Snell's law (\ref{eq:Snells_Law}), the GEOFA system reduces to one equation for the last remaining variable $x_n$, reading
\[
    c_{n+1} \; \frac{x_{n} - x_{n-1} }{ \sqrt{(x_{n} - x_{n-1})^2 + (d_n - z_{n-1})^2}} = c_n \; \frac{x_{n+1} - x_{n} }{\sqrt{(x_{n+1} - x_{n})^2 + (z_{n+1} - d_n)^2} } \, .        
\]

\paragraph{Elliptic Medium Boundary}

For an elliptic medium boundary $B_n$, the medium boundary function takes the form $b_n(x) = \pm \, b \, \sqrt{1 - x^2 / a^2 }$, where the sign depends on the outward or inward curvature of the medium boundary. Therefore, the equations (\ref{eq:medium_boundary_equation}) and (\ref{eq:tangent_of_medium_boundary_equation}) turn into
    \begin{align}
        z_n = b_n(x_n) = \pm \, b \, \sqrt{1 - x_n^2 / a^2 } \, , \label{elliptic_medium_boundary_equation}\\
        \tan(\alpha_n) = b_n'(x_n) = \frac{\mp \, b x_n/a^2}{\sqrt{1 - x_n^2 / a^2 }}\, .  \label{elliptic_medium_boundary_slope}
    \end{align}

\section{Numerical and Experimental Results} \label{sect_4:numerical_simulations}

This section presents numerical simulations and phantom-based experiments evaluating the GEOFA algorithm (\cref{Alg:Detailed_multilayered_focusing_Algorithm}) as part of our ToF correction approach, see \cref{eq:acoustic_lens_receive}. \cref{sect_4_1_ToF_measurements} reports ToF and focusing delay simulations performed with the acoustic toolbox k-Wave \cite{First_kWave_Paper}. In \cref{sect_4_2:ToF_errors_depending_on_med_bound_points}, GEOFA is compared against a similar method from \cite{MINEO2022106747}. Finally, \cref{sect_4_3_fat_layer_image_improvement} demonstrates ultrasound image aberration correction in a phantom-based experimental setup, where an aberrating medium is interposed between the ultrasound transducer and the phantom.

\subsection{Time of Flight and Focusing Delay Simulations} \label{sect_4_1_ToF_measurements}

ToFs and focusing delays calculated using GEOFA, SRT and HMFA are compared to 2D simulations conducted using the acoustic toolbox k-Wave.

\begin{figure}[ht!]
    \centering
    \includegraphics[width=\textwidth]{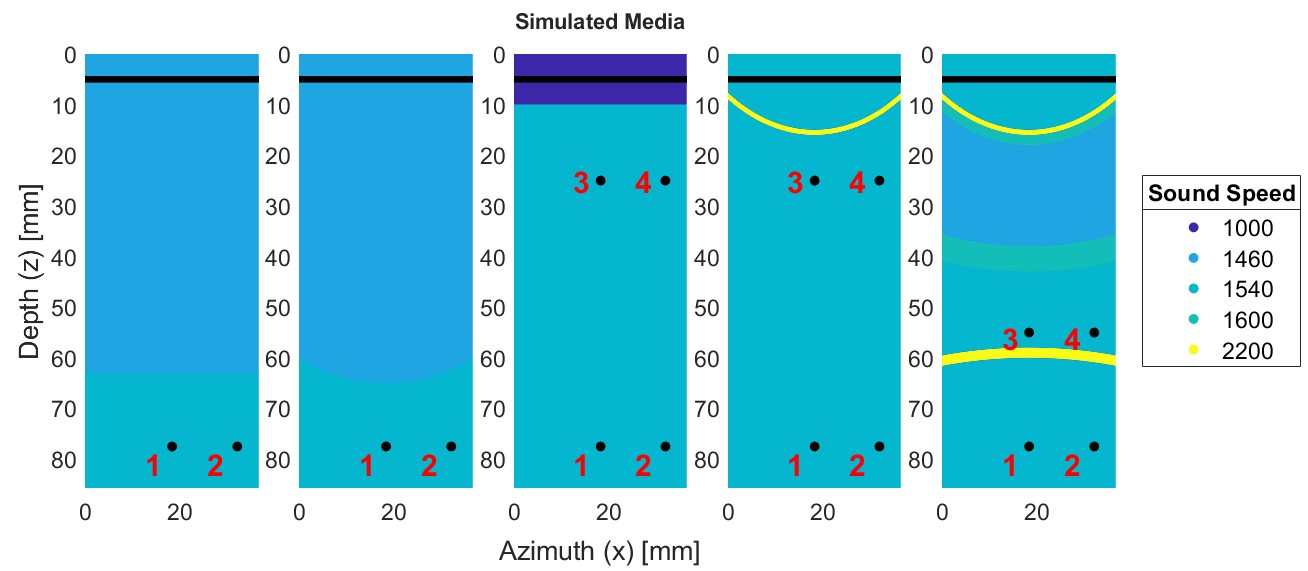}
    \caption{Layered media simulation settings with sources (black dots) of spherical waves inside the medium and sensors (black line) in the first layer.}
    \label{fig:ToF_simulations_Settings}
\end{figure}

\begin{table}[htbp]
\centering
\caption{Sound speeds of simulated tissues.}
\label{tab:tissue_sound_speed}
\begin{tabular}{lcccccc}
\toprule
Tissue & \makecell{Subcutaneous\\fat} & \makecell{Maternal\\skin} & \makecell{Rectus\\abdominis\\(muscle)} & \makecell{Amniotic\\fluid} & \makecell{Fetal\\calvaria\\(skull)} & \makecell{Fetal\\brain} \\
\midrule
$c$ [\si{m/s}] & 1460 & 1600 & 1600 & 1540 & 2200 & 1540 \\
\bottomrule
\end{tabular}
\end{table}

\paragraph{Simulation settings} \label{sect_4_1_1_simulation_settings}

The five settings used in the k-Wave simulations are visualized in \cref{fig:ToF_simulations_Settings}. In each simulation, an acoustic wave is emitted from one of the sources (black dots in \cref{fig:ToF_simulations_Settings}). The acoustic wave is detected by sensors distributed along the top of the
medium (black line in \cref{fig:ToF_simulations_Settings}). By time-reversal invariance \cite{Cass_Wu_Fink_closed_time_reversal_cavities}, the ToF for a wave traveling through the medium equals the ToF of a wave on the reverse path. This approach is advantageous as it directly provides the ToFs from one point inside the medium (i.e., the focus point) to multiple points in the first layer (i.e., the transducer elements), which is the typical configuration in ultrasound imaging. The focusing delays (\ref{eq:acoustic_lens}) required for focusing on the origin of the acoustic wave are then derived from the ToF simulation results.

To ensure our simulation accurately reflects practical applications, we choose sound speeds and geometries occurring in real medical ultrasound \cite{Tissue_Properties_3,Tissue_Properties_1,Tissue_Properties_2}, see \cref{tab:tissue_sound_speed}. The simulated media are visualized in \cref{fig:ToF_simulations_Settings}. The first setting consists of a fat layer on top of a tissue layer, separated by a flat medium boundary. The second setting is similar, differing only by an elliptic medium boundary. In the third and fourth setting, we emulate a flat and an elliptic transducer cover. Finally, the fifth setting resembles fetal brain imaging, comprising a realistic sequence of tissue layers (maternal skin, subcutaneous fat, abdominal muscle, amniotic fluid, fetal skull, and brain) beneath an elliptic transducer cover. Details regarding the simulation geometry and the k-Wave simulation settings are provided in \cref{sect:B_kWave_Simulation_Details}. 

\paragraph{Results: Time of Flight and focusing delay accuracy} 
\label{sect_4_1_2_ToF_Results}

In the following, ToFs and focusing delays calculated using HMFA, SRT, and GEOFA are compared to k-Wave simulation results.

\begin{definition}[ToF and focusing delay errors]
    Assume that ToFs and focusing delays (\ref{eq:acoustic_lens}) obtained by simulations discussed in this section are given. Then, for any of the methods HMFA, SRT, or GEOFA, the ToF error is defined as the difference between the method's ToFs and the simulated ToFs. Analogously, the focusing delay error is defined as the difference between the focusing delays calculated by any of the above methods and the focusing delays obtained by the simulation results.
\end{definition}

\begin{figure}[ht!]
    \centering
    
    \begin{minipage}{.485\textwidth}
        \includegraphics[width=\textwidth]{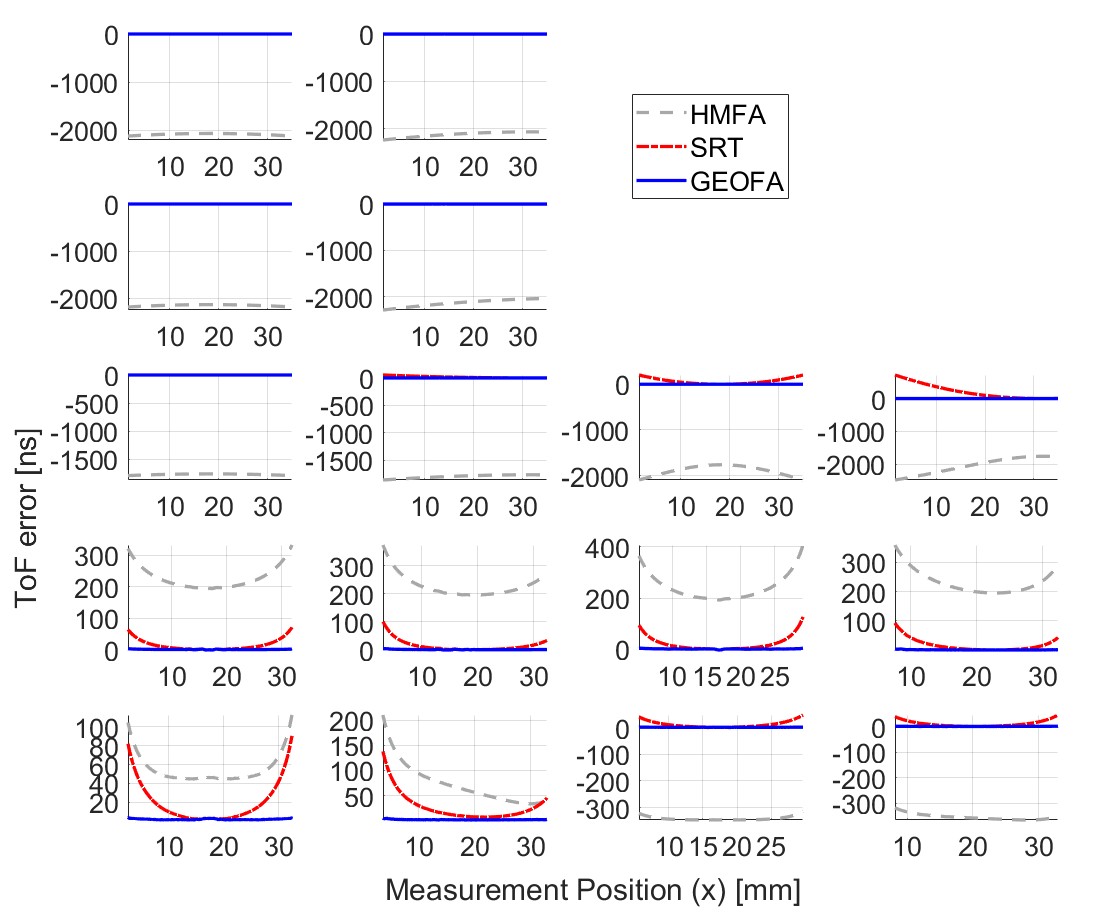}
    \end{minipage}
    \hfill
    \begin{minipage}{.485\textwidth}
        \includegraphics[width=\textwidth]{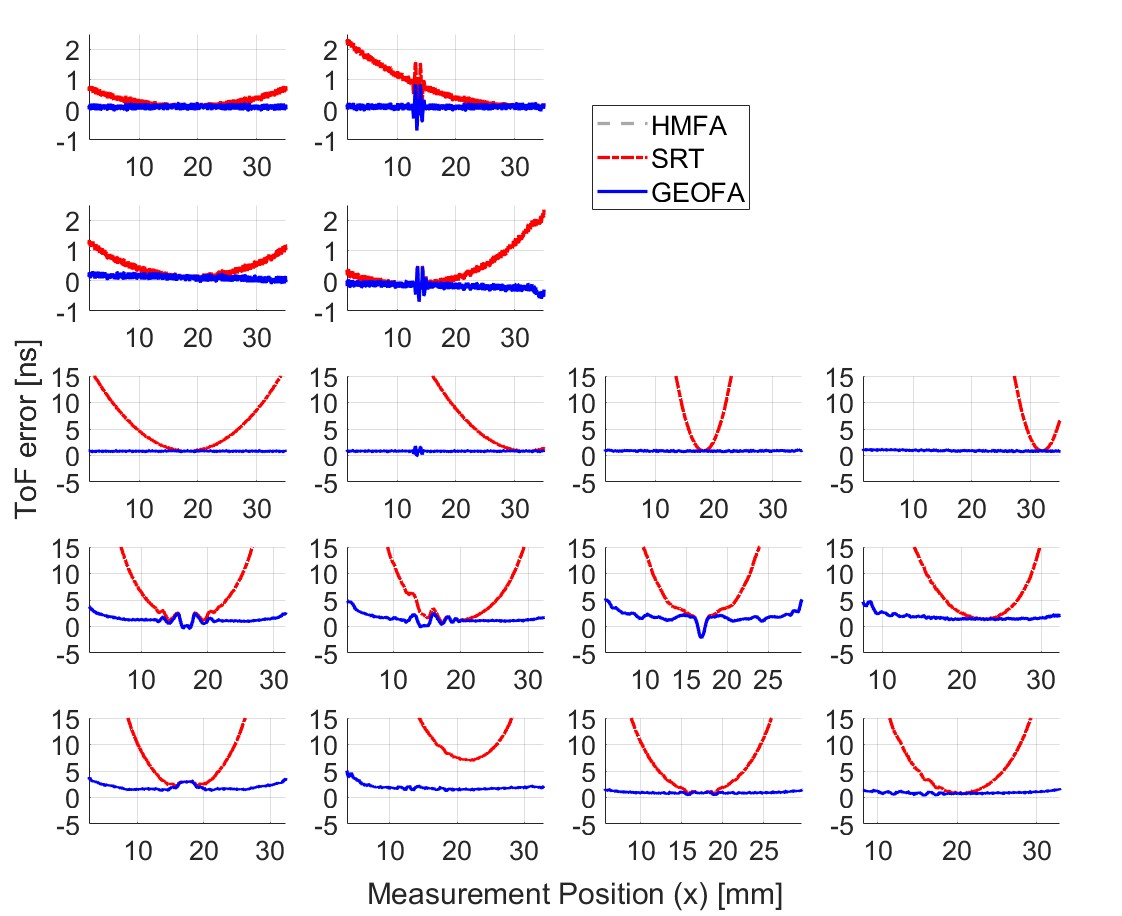}
    \end{minipage}
    \caption{ToF errors of HMFA, SRT, and GEOFA for the five simulation settings of \cref{fig:ToF_simulations_Settings}. Each panel shows the error as a function of the lateral sensor position. Within each block, the rows correspond, from top to bottom, to Settings $1$ to $5$, and the two columns correspond to the two source positions, ordered according to the source numbers in \cref{fig:ToF_simulations_Settings}. Left block: full error range. Right block: identical data with a re-scaled $y$-axis for a detailed assessment of the SRT and GEOFA errors. The small-scale oscillations of the SRT and GEOFA curves reflect oscillations in the simulated ground truth, which are discussed in the text.}
    \label{fig:ToF_errors}
\end{figure}

\begin{figure}[ht!]
    \centering
    
    \begin{minipage}{.485\textwidth}
        \includegraphics[width=\textwidth]{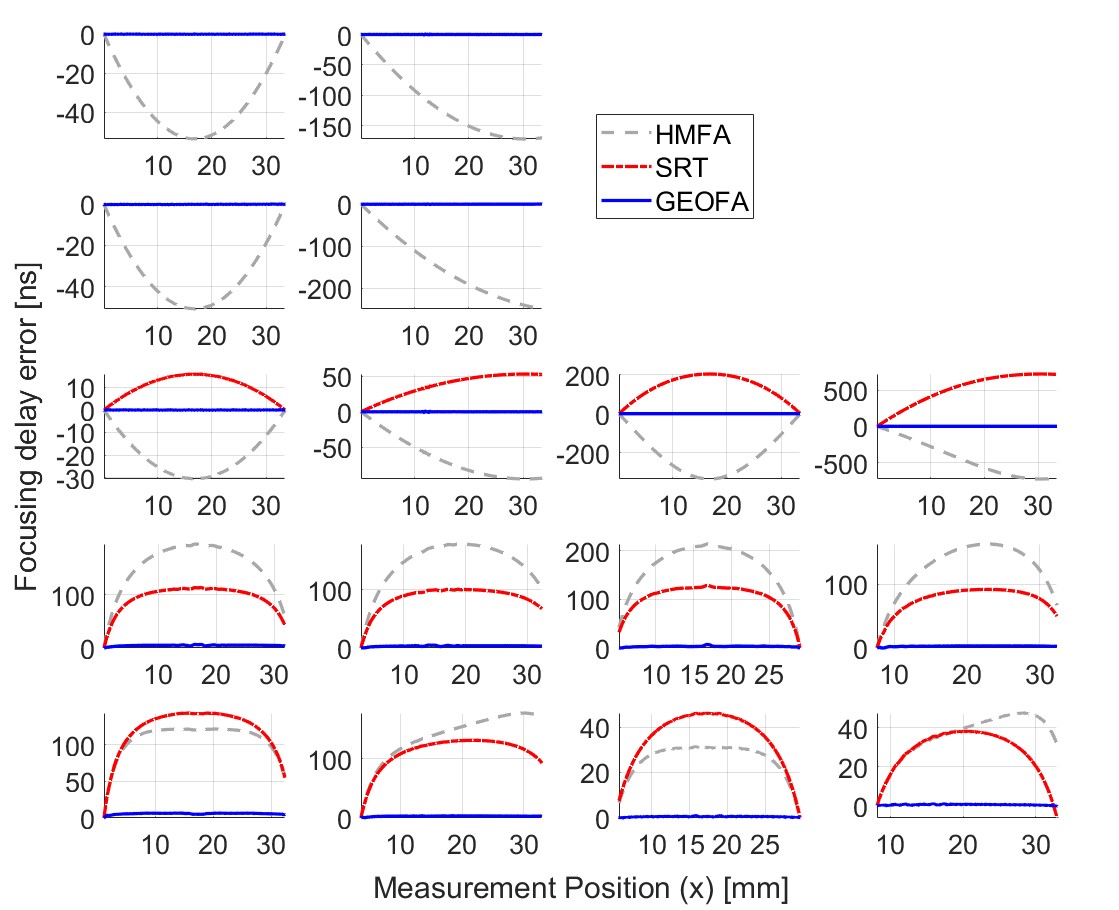}
    \end{minipage}
    \hfill
    \begin{minipage}{.485\textwidth}
        \includegraphics[width=\textwidth]{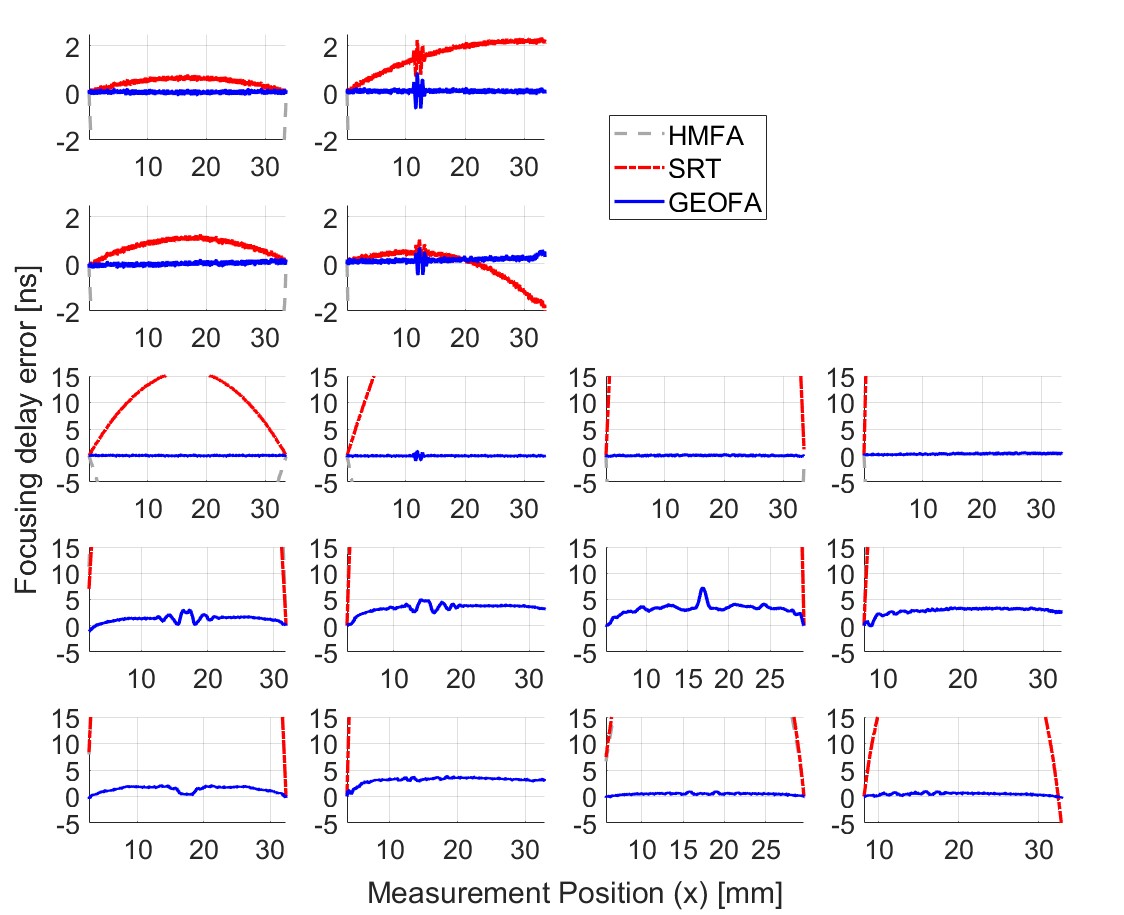}
    \end{minipage}
    \caption{Focusing Delay errors of HMFA, SRT, and GEOFA for the five simulation settings of \cref{fig:ToF_simulations_Settings}. Each panel shows the error as a function of the lateral sensor position. Within each block, the rows correspond, from top to bottom, to Settings $1$ to $5$, and the two columns correspond to the two source positions, ordered according to the source numbers in \cref{fig:ToF_simulations_Settings}. Left block: full error range. Right block: identical data with a re-scaled $y$-axis for a detailed assessment of the SRT and GEOFA errors. The small-scale oscillations of the SRT and GEOFA curves reflect oscillations in the simulated ground truth, which are discussed in the text.}
    \label{fig:focusing_delay_errors}
\end{figure}

\begin{remark}[Significance threshold for focusing delay errors]
\label{rem:threshold}
    In digital beamformer design, tolerable focusing delay errors between $1/16$ and $1/32$ of the fundamental period (equivalently, of the wavelength) are commonly reported in the literature \cite{Fritsch2006, Steinberg1992}, with $1/32$ stated as the desired accuracy for high image quality \cite{Alexandru2010, Agarwal2016}. To be on the safe side, we adopt the strict end of this range and consider focusing delay errors exceeding $1/32$ of the fundamental period as significant, which at $5\,\mathrm{MHz}$ corresponds to $6.25\,\mathrm{ns}$. For quantitative analyses of the effect of delay errors on the resulting image quality, we refer to \cite{Beaver1977, Magnin1981, Focusing_error_source}.
\end{remark}

\cref{fig:ToF_errors} illustrates the ToF errors, while the focusing delays are visualized in \cref{fig:focusing_delay_errors}. In the right panel of both figures, the $y$-axis is re-scaled for a more detailed assessment of the GEOFA errors.

The HMFA ToF errors in Settings $1$ and $2$ are similar, since the small sound speed contrast between fat and tissue makes refraction negligible. In Setting $4$, the sign of the error reverses, since all sound speeds exceed the $1540\,\text{m/s}$ assumed in HMFA, whereas they are below this value in the first three settings. The thinner aberrating layer offsets the larger sound speed difference, which leads to a lower absolute error. For the right-shifted sources, the ToF error increases due to a longer acoustic path through the aberrating media. Additionally, in Settings $3$, $4$, and $5$, a larger angle of incidence with respect to the medium boundary normal amplifies refraction, see~\cref{eq:Snells_Law}. Regarding focusing delay errors of HMFA, Settings $1$ and $2$ are nearly identical. Despite smaller absolute ToF errors, Settings $3$, $4$, and $5$ exhibit substantially larger focusing delay errors than Settings $1$ and $2$, because the ToF error curve is steeper. Since focusing delays are defined as differences relative to the maximum ToF~(see~Definition~\ref{def:Transmit_Focusing}), the curvature of the ToF error profile matters more than its magnitude. In all settings, focusing delay errors are significant, exceeding the threshold by $1$ to $2$ orders of magnitude, which corresponds to $\frac{1}{3}$ to $3$ wavelengths at $5\,\mathrm{MHz}$.

SRT reduces the HMFA ToF error substantially in most settings. In Settings $1$ and $2$, the ToF errors of SRT are small enough such that the resulting focusing delay errors remain insignificant. However, in Settings $3$, $4$, and $5$, the neglect of Snell's law of refraction in SRT (\ref{eq:Snells_Law}) results in a significant ToF and focusing delay error with a characteristic parabolic profile, which is amplified for right-shifted sources and, in Settings $3$ and $4$, near-field positions. The focusing delay errors in Settings $3$, $4$, and $5$ exceed the significance threshold by approximately one order of magnitude, which corresponds to $\frac{1}{3}$ wavelength at $5\,\mathrm{MHz}$. Consequently, for large sound speed variations, such as those introduced by bone or transducer cover materials, SRT produces significant focusing delay errors, which can only be corrected by accounting for Snell's law of refraction, as GEOFA does.

In Settings $1$ to $3$ and $5$, the ToF and focusing delay errors of GEOFA are entirely below the significance threshold of $6.25$~ns. Only in Setting $4$ does the error slightly exceed this rather strict threshold at times, although it remains within the same order of magnitude and below $\frac{1}{24}$ of a wavelength at $5\,\mathrm{MHz}$. The errors in Settings $4$ and $5$ are presumably caused by wave-like phenomena at the transducer cover, such as diffraction or interference, which cannot be modeled using geometrical acoustics, see \cref{sect:2_3_geometrical_acoustics}. It is important to note that the simulation in Setting $5$ was performed with a finer grid for an improved precision required due to the multiple curved medium boundaries in the simulation, see \cref{sect:B_kWave_Simulation_Details}.

The small-scale oscillations visible in the SRT and GEOFA curves on the re-scaled panels of \cref{fig:ToF_errors,fig:focusing_delay_errors} are not properties of the focusing algorithms but reflect fluctuations in the simulated ground truth. Since the ToFs calculated by HMFA, SRT, and GEOFA are smooth, non-oscillating functions of the sensor position, any non-smooth, oscillatory component of the error curves must originate from the simulated ToF, which is confirmed by the same oscillation pattern appearing in the error curves of all three methods. These oscillations are highly localized and remain within a few nanoseconds, which suggests that they are caused by wave-like phenomena in the simulation, such as diffraction or interference, which cannot be modeled using geometrical acoustics, see \cref{sect:2_3_geometrical_acoustics}. We observe that flat medium boundaries only occasionally cause such small oscillations without any systematic deviation from GEOFA, whereas curved boundaries, in particular the two-layered curved transducer cover, produce both stronger oscillations and a small systematic deviation that grows toward the more strongly curved part of the cover. A minor additional contribution stems from temporal quantization of the underlying simulation. The time discretization used is $dt = 0.25$\,ns, see \cref{sect:B_kWave_Simulation_Details}. This should be kept in mind when discussing artifacts on the order of a few time discretization steps.

To summarize, these k-Wave simulation results show that using GEOFA, ToFs can be calculated highly accurately, and focusing delay errors of HMFA and SRT can be reduced to near insignificance, see Remark~\ref{rem:threshold}. This result strongly supports our approach of using GEOFA for focusing delay error correction.

\subsection{Comparison of GEOFA with State-of-the-Art Geometrical 
Acoustics Method} \label{sect_4_2:ToF_errors_depending_on_med_bound_points}

This section compares GEOFA with the algorithm introduced in~\cite{MINEO2022106747}, which we refer to as MINEO. MINEO is a state-of-the-art geometrical acoustics method for ToF calculation in known layered media. Numerical calculations are used to compare GEOFA and MINEO, adding to the theoretical comparison in \cref{sect:3_4_MAFA}.

Each algorithm’s nonlinear system of equations is solved using the trust-region-dogleg method, implemented in MATLAB's \texttt{fsolve}. No limit on the number of iterations was imposed, but the solver terminates when no further progress is achieved. For the initial guess, we use the medium boundary points obtained by SRT. Since both algorithms use the same iterative solver and the same initial points, their intermediary iterates and final results are directly comparable. In all cases where both algorithms converged, the resulting ToFs agreed up to machine precision, multiple orders of magnitude below the significance level of $6.25$\,ns defined in Remark~\ref{rem:threshold}.

\begin{figure}[ht!]
    \centering
    \begin{minipage}{.425\textwidth}
        \includegraphics[width=\textwidth]{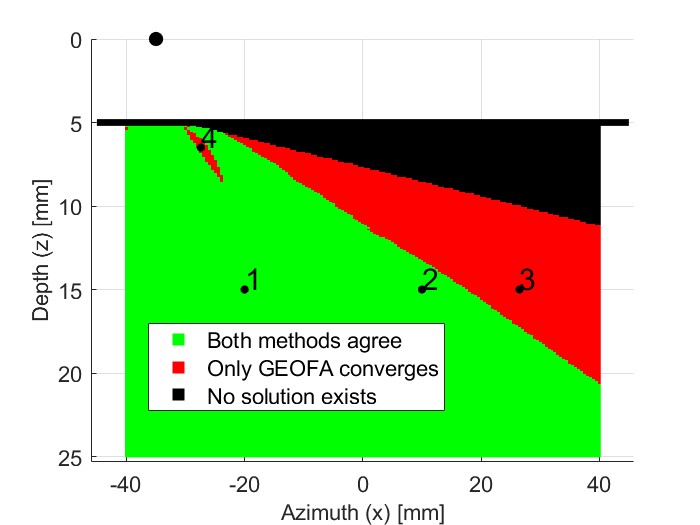}
    \end{minipage}
    \hfill
    \begin{minipage}{.525\textwidth}
        \includegraphics[width=\textwidth]{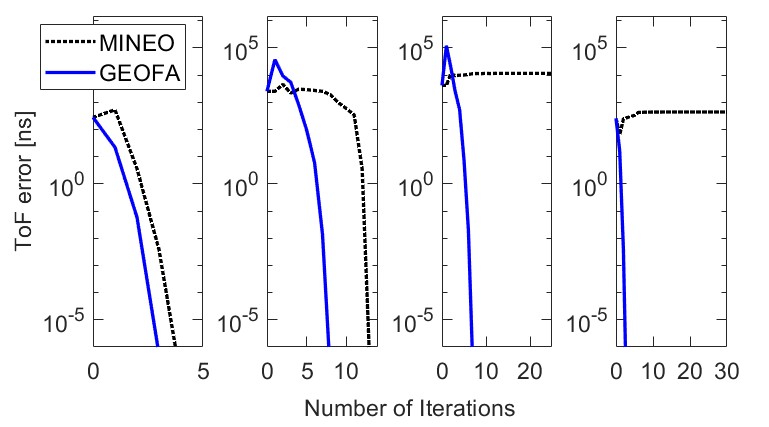}
    \end{minipage}
    \caption{Straight transducer cover setting with source at $(-35,0)$\,mm. Left: convergence map. Right: ToF error per iteration for target points $1$ to $4$.}
    \label{fig:Algorithm_comparison_1}
\end{figure}

\begin{figure}[ht!]
    \centering
    \begin{minipage}{.425\textwidth}
        \includegraphics[width=\textwidth]{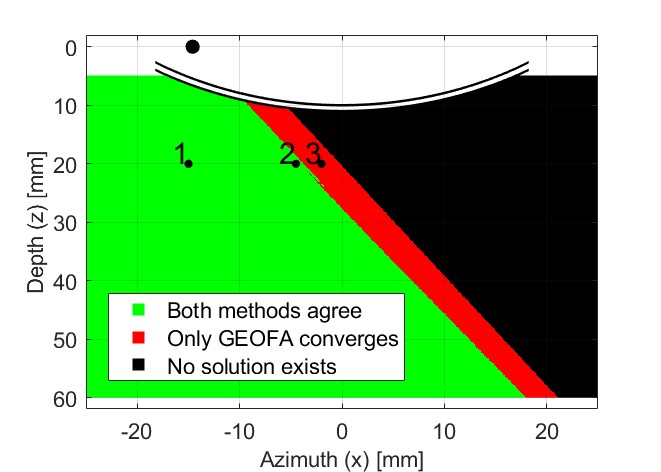}
    \end{minipage}
    \hfill
    \begin{minipage}{.525\textwidth}
        \includegraphics[width=\textwidth]{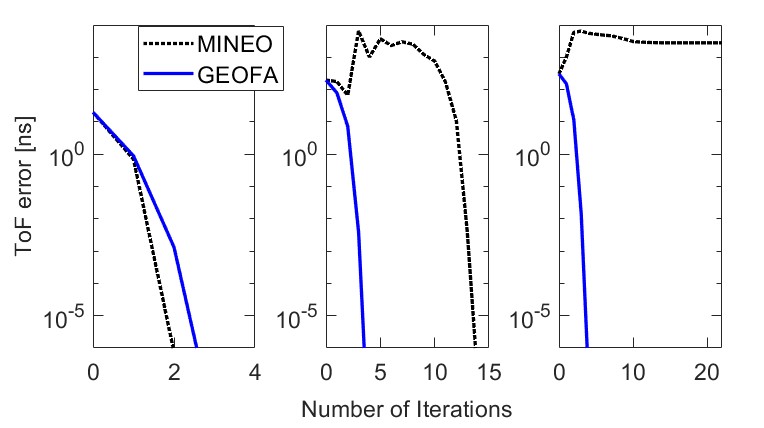}
    \end{minipage}
    \caption{Elliptic transducer cover setting with source at $(-14.6,0)$\,mm. Left: convergence map. Right: ToF error per iteration for target points $1$ to $3$.}
    \label{fig:Algorithm_comparison_2}
\end{figure}

\begin{figure}[ht!]
    \centering
    \begin{minipage}{.425\textwidth}
        \includegraphics[width=\textwidth]{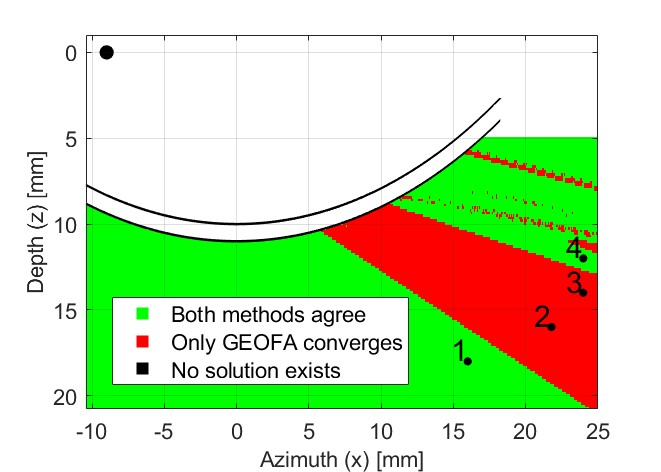}
    \end{minipage}
    \hfill
    \begin{minipage}{.525\textwidth}
        \includegraphics[width=\textwidth]{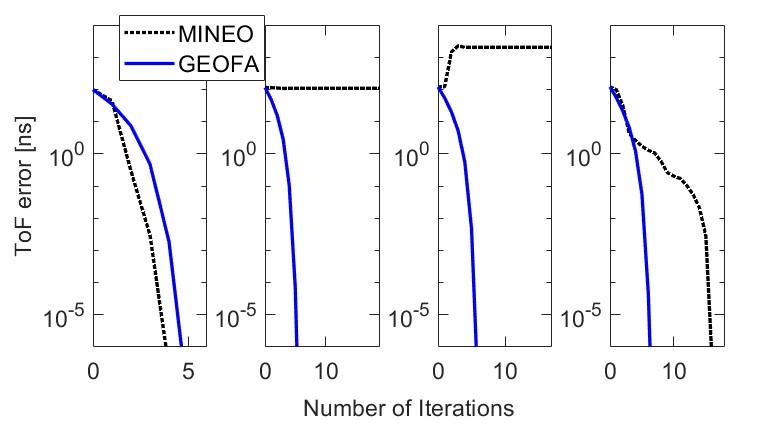}
    \end{minipage}
    \caption{Elliptic transducer cover setting with source at $(-9,0)$\,mm. Left: convergence map. Right: ToF error per iteration for target points $1$ to $4$.}
    \label{fig:Algorithm_comparison_3}
\end{figure}

We evaluate both methods on Settings~$3$ and~$4$ from~\cref{sect_4_1_ToF_measurements}, visualized in~\cref{fig:ToF_simulations_Settings}. For each setting, a source point is fixed in the first layer (indicated by the black dot in the left panels of Figures~\ref{fig:Algorithm_comparison_1} to~\ref{fig:Algorithm_comparison_3}), and target points are chosen on a finely discretized grid in the last layer. The color at each target point indicates the convergence outcome: green if both methods converge to the same solution, red if only GEOFA converges, and black if no solution exists due to total reflection. The numbered dots mark selected target points whose convergence behavior is shown in the right panels, which plot the ToF error after each iteration.

The straight medium boundary of the first setting (\cref{fig:Algorithm_comparison_1}) makes the geometry translation invariant along the $x$-axis. Moreover, the result is symmetric about the vertical line through the source. Therefore, a single source position suffices to obtain the general error pattern. In this setting, there exists a triangular region of target points adjacent to the total reflection border, where MINEO fails to converge but GEOFA succeeds. In addition, there is a small area closer to the medium boundary where GEOFA converges but MINEO does not.

The second and third settings (\cref{fig:Algorithm_comparison_2,fig:Algorithm_comparison_3}) use an elliptical transducer cover. For a source near the edge of the transducer (\cref{fig:Algorithm_comparison_2}), a large region is inaccessible due to total reflection. In a slice adjacent to the total reflection border, MINEO fails to converge, while GEOFA still succeeds. As the source moves toward the center (\cref{fig:Algorithm_comparison_3}), the total-reflection region shrinks and, in \cref{fig:Algorithm_comparison_3}, has vanished entirely, but a region of non-convergence of MINEO persists.

The convergence plots in the right panels of Figures~\ref{fig:Algorithm_comparison_1} to~\ref{fig:Algorithm_comparison_3} illustrate this in more detail. For target points well within the green region (e.g., point~$1$ in all figures), both methods converge within the first five iterations. As the target point approaches the total-reflection border, GEOFA continues to converge, although requiring more iterations, while MINEO stagnates at a large ToF error. This behavior is consistent across all tested configurations: GEOFA converges up to the physical limit imposed by total reflection, whereas MINEO might fail in a transitional region of strong refraction. This observation can be explained by analyzing the GEOFA and the MINEO systems of equations.

The reason for the non-convergence of MINEO lies in the functional form of its equations rather than in the choice of parametrization, which for a horizontal boundary coincides with ours. In the MINEO system, Snell's law is enforced by computing the post-refraction angle through the inverse sine, which is real-valued only for sub-critical incidence angles and whose derivative diverges as the critical angle of total reflection is approached. Moreover, in each layer, the lateral displacement of the ray is proportional to $\tan \theta$, with $\theta$ being the angle of the ray with respect to the layer normal. Close to the critical angle, refraction bends the ray toward shallower angles with the interface, the refracted angle $\theta'$ approaches $90^\circ$, and the lateral displacement grows without an upper bound. The MINEO system therefore contains a singularity exactly at the physical limit of total reflection, and whenever the iterates of the Newton-type method approach or cross the critical angle, its Jacobian degenerates, the equations are no longer real-valued, and the solver stagnates. This explains both regions of non-convergence in \cref{fig:Algorithm_comparison_1}. In the region adjacent to the total reflection border, the solution itself lies close to the critical angle. In the small region close to the medium boundary, the initialization is super-critical although the solution is not. For target point $4$, the straight ray to the SRT medium boundary point, which both methods use as initialization, has a super-critical incidence angle, while the incidence angle of the solution is sub-critical. Since for targets located closely below the medium boundary the updates of the transmitted angle additionally scale inversely with the distance of the target below the boundary, the iterates are driven across the singularity and the solver fails. The GEOFA system does not contain this singularity, since the sines of the angles are themselves unknowns, Snell's law \eqref{eq:Snells_Law} is a linear equation in these unknowns, and the remaining equations are normalized inner products, which are bounded together with their derivatives for arbitrary iterates. Therefore, GEOFA remains real-valued and well conditioned up to the physical limit of total reflection.

The authors of~\cite{MINEO2022106747} suggest using bisection to improve the initial guess before applying a Newton-type solver. This addresses convergence issues in industrial material testing, where multiple layers with large sound speed differences can occur~\cite{MINEO2022106747}. In medical ultrasound imaging, however, the sound speeds of most tissue types are close to the $1540$\,m/s assumed by HMFA, except for transducer cover materials and bone, which have a substantially different sound speed. Consequently, the number of strongly refracting interfaces is small, and the solution of the nonlinear system can be found using SRT as an initial guess. The results presented here confirm that under these conditions, GEOFA converges reliably with SRT as initial guess alone, without requiring preliminary bisection iterations.

\subsection{Experimental Evaluation of ToF-Based Aberration Correction}\leavevmode\par
\label{sect_4_3_fat_layer_image_improvement}
This section evaluates the performance of our image reconstruction approach, which uses GEOFA for receive focusing in the presence of layered medium aberrations.

\paragraph{The Experimental Setup}

In our experiment, we use the Linear Array Ultrasound Probe 9L-D (GE HealthCare) with a transmit center frequency of $6.75$~MHz, a focal depth of $55$~mm, and a transmit aperture size of $12.9$~mm. Data were acquired using a commercial GE system, which does not allow access to the exact transmit waveform or pulse sequence. We use a Sono404 (Sun Nuclear) phantom with a sound speed of $1540$~m/s. Furthermore, two 3D-printed aberrators of silicone 40A printed with a Form 3B printer by Formlabs are used. The first aberrator is a $3$~mm thick cuboid, the second a $0$ to $3$~mm thick wedge.

The sound speed of the aberrator material, $c_{\text{aberrator}}$, was previously unknown and reconstructed in a first experiment. For this purpose, the cuboid aberrator is used to record an aberrated image. Subsequently, a ground truth image is recorded by carefully removing the aberrator without disturbing the setup. The space previously occupied by the aberrator is filled with $4.5$\% saline solution, which has a sound speed of approximately $1540$ m/s at room temperature \cite{4_5_per_cent_Salt_water_at_room_temperature_1540_Soundspeed}. Both images are reconstructed using HMFA beamforming, implemented in the Ultrasound Toolbox (USTB) \cite{USTB}. The cuboid is chosen for reconstruction because its uniform thickness yields a simple, experimentally stable relationship between the depth shift in the image and $c_{\text{aberrator}}$. The wedge, with its spatially varying thickness, provides a more challenging geometry for aberration correction. An aberrated image is recorded, and a ground truth image is obtained after replacing the aberrator with $4.5$\% saline solution. The sound speed of the aberrator material, $c_{\text{aberrator}} = 1030.2 \pm 5.4$\,m/s (mean $\pm$ standard deviation), was reconstructed by comparing the HMFA images with and without cuboid, see \cref{sect:C_SoS_reconstruction}.

\paragraph{Quantitative evaluation of reconstructed images using lateral wire target profiles}

\begin{figure}[ht!]
    \centering
    \includegraphics[width = \textwidth]{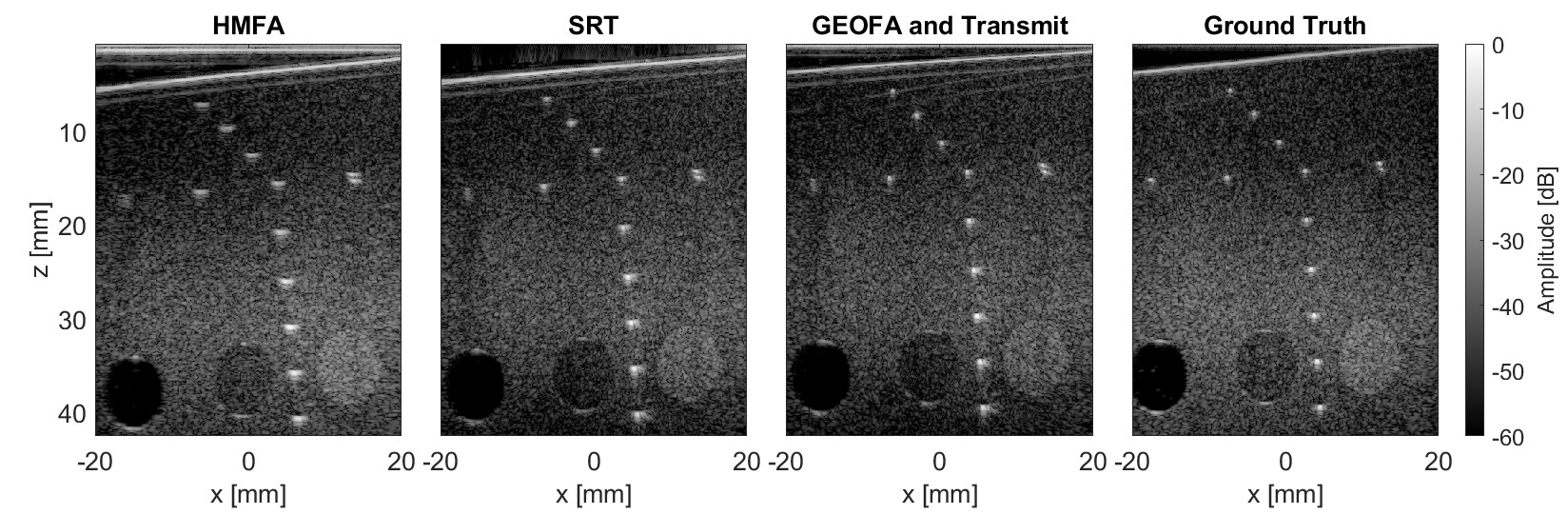}
    \caption{Reconstructed Ultrasound Images of the wedge aberrator (first three) and the ground truth (fourth image). Left: HMFA, Center Left: SRT correction, Center Right: GEOFA and Transmit correction, Right: HMFA without aberrator, for comparison. The bright horizontal bands in the first $15$\,mm are reverberation artifacts caused by multiple reflections between the medium boundaries and the transducer, see the discussion at the end of this section. The doubled appearance of the rightmost wire target in the horizontal row is caused by two closely spaced scatterers in the phantom.}
    \label{fig:9L_Fat_layer_images}
\end{figure}

\begin{figure}[ht!]
    \centering
    \begin{minipage}{.315\textwidth}
        \includegraphics[width=\textwidth]{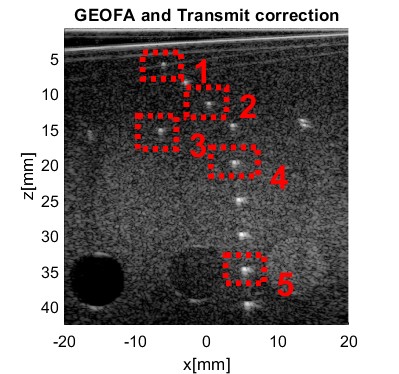}
        \includegraphics[width =\textwidth]{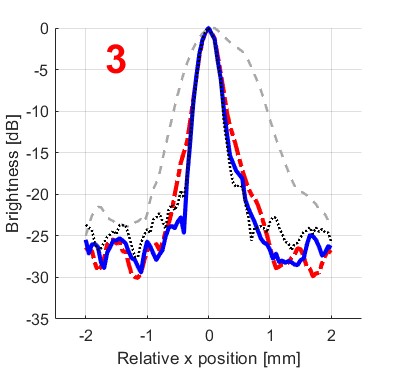}
    \end{minipage}
    \hfill
    \begin{minipage}{.665\textwidth}
        \centering
        \begin{minipage}{.485\textwidth}
            \includegraphics[width=\textwidth]{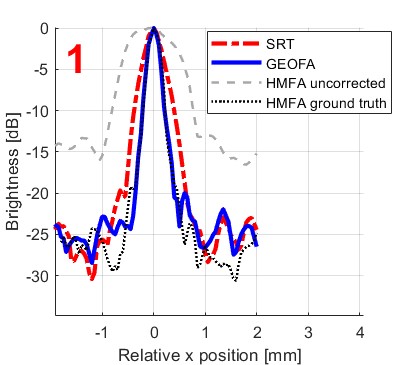}
        \end{minipage}
        \hfill
        \begin{minipage}{.485\textwidth}
            \includegraphics[width=\textwidth]{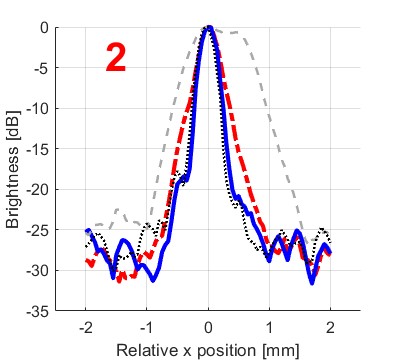}
        \end{minipage}\centering
        \hfill
        \begin{minipage}{.485\textwidth}
            \includegraphics[width=\textwidth]{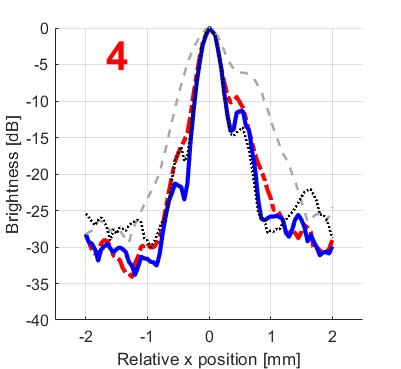}
        \end{minipage}
        \hfill
        \begin{minipage}{.485\textwidth}
            \includegraphics[width=\textwidth]{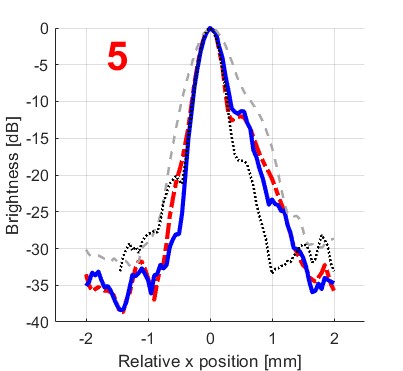}
        \end{minipage}
    \end{minipage}
    \caption{Wire target profiles evaluated at multiple wire targets.}
    \label{fig:Beam_Profile_evaluation_using_Wire_Targets}
\end{figure}

Image quality is assessed by analyzing the wire targets visible in \cref{fig:9L_Fat_layer_images}. The leftmost panel, showing the uncorrected HMFA image with the wedge aberrator present, exhibits clear lateral blurring. To quantify this effect, lateral wire target profiles are evaluated as follows: a rectangular region is selected around each wire target, as illustrated in the top left panel of \cref{fig:Beam_Profile_evaluation_using_Wire_Targets}, and the lateral wire target profile along the azimuth (x) direction is obtained by computing the maximum intensity projection along the depth (z) direction. The resulting lateral wire target profiles for all evaluated wire targets are used as the local image quality measure.

\paragraph{Our ToF correction approach}

To better approximate the underlying ToF in the focusing delays used for receive focusing (see \cref{eq:acoustic_lens_receive}), we apply corrections to both the transmit and receive components.

\begin{remark}[GEOFA and Transmit ToF correction] \label{rem:transmit_and_receive_correction}
GEOFA is used to correct the receive ToF. The transmit ToF correction is computed, for each receive focus point, along the vertical ray through that point, using $c_{\text{aberrator}}$ in the aberrating layer instead of the default $c = 1540\,\text{m/s}$ assumed in HMFA beamforming. In contrast, the SRT comparison method is applied, as in the state-of-the-art literature \cite{jaeger2015computed,sanabria2018spatial,vraalstad2024coherence}, to the receive ToF only and does not contain a transmit correction.
\end{remark}

\paragraph{Comparison to beamforming using a constant sound speed}

\begin{figure}
    \centering
    \includegraphics[width=.6\linewidth]{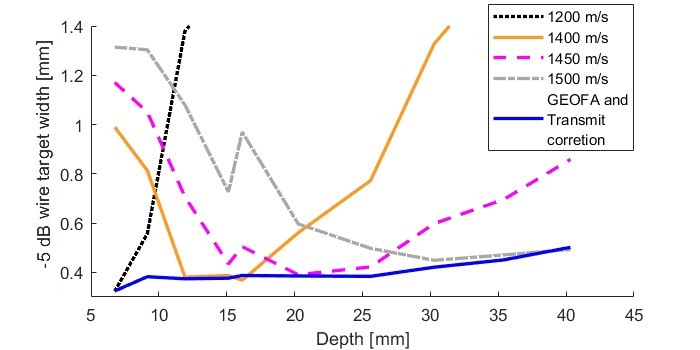}
    \caption{$-5$ dB lateral wire target widths for HMFA images reconstructed with a range of different constant sound speeds, compared to our GEOFA and Transmit correction, which achieves minimal wire target width across all depths.}
    \label{fig:3D_print_material_SoS_measurement}
\end{figure}

In this paragraph, we investigate whether a globally optimized constant sound speed, potentially reconstructable with previously proposed methods \cite{SoS_estimate_steering}, can improve HMFA image quality in the wedge-aberrated setting. For this, HMFA images are reconstructed using a range of different constant sound speeds. The resulting $-5$~dB lateral wire target widths of all wire targets are plotted in \cref{fig:3D_print_material_SoS_measurement}. These results are compared with our ToF correction approach, see Remark \ref{rem:transmit_and_receive_correction}.

The wire target widths reveal a clear depth-dependent behavior: each constant sound speed produces narrow wire targets at some depths while significantly degrading resolution at others. Sound speeds closer to $c_{\text{aberrator}}$ yield better resolution at shallow depths, while values closer to $1540$ m/s perform better at greater depths, consistent with the diminishing influence of the aberrator. Only our ToF correction approach achieves consistently narrow wire target profiles across all depths.

\paragraph{Comparison of GEOFA, SRT, and HMFA}

We apply two reconstruction methods, both of which we implemented in the USTB, and compare the resulting images with the wedge-aberrated and ground-truth HMFA images. First, as a baseline comparison with a state-of-the-art method, SRT is applied to the receive ToF part of the receive focusing delays (\ref{eq:acoustic_lens_receive}). Second, GEOFA is applied to the receive ToF and a transmit ToF correction is additionally considered, see Remark \ref{rem:transmit_and_receive_correction}.

The resulting images are shown in \cref{fig:9L_Fat_layer_images}. A lateral blurring of the wire targets can be observed in the aberrated HMFA image, which is already significantly improved by SRT. The combination of GEOFA and Transmit ToF correction improves upon SRT, especially in the near field, and only minor aberrations remain compared to the ground truth image.

\begin{figure}[ht!]
    \centering
    \begin{minipage}{.475\textwidth}
        \includegraphics[width=\textwidth]{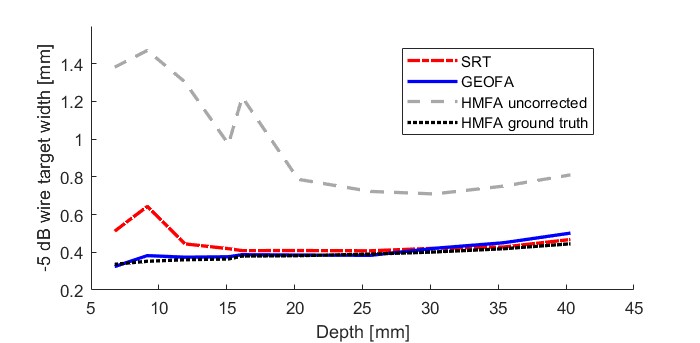}
    \end{minipage}
    \hfill
    \begin{minipage}{.475\textwidth}
        \centering
        \includegraphics[width=\textwidth]{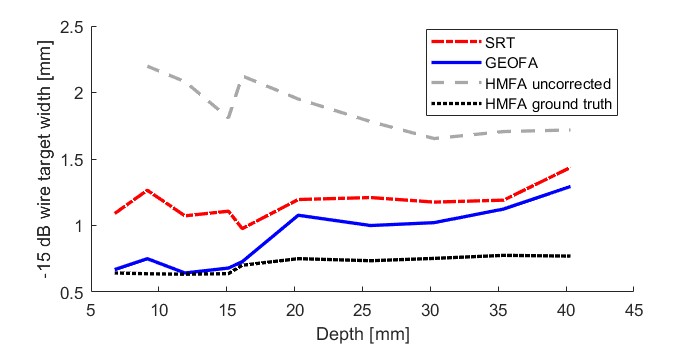}
    \end{minipage}
    \caption{Lateral wire target profile widths for varying depth at relative brightness of $-5$ dB (left) and $-15$ dB (right).}
    \label{fig:Beam_Width_over_depth}
\end{figure}

For a quantitative evaluation, the wire target profiles are shown in \cref{fig:Beam_Profile_evaluation_using_Wire_Targets}. SRT yields significantly narrower wire target profiles than HMFA. Our GEOFA and transmit correction further improves upon SRT, approaching the ground truth for targets $1$ through~$4$. For target~$5$, a deviation from the ground truth is observed, which we attribute to retrospective transmit beamforming (RTB), discussed in more detail below.

\cref{fig:Beam_Width_over_depth} shows a quantitative wire target width evaluation at $-5$\,dB and $-15$\,dB for all wire targets, except the leftmost and rightmost target in the horizontal row, which are excluded for the same reasons as in the sound speed reconstruction, see \cref{sect:C_SoS_reconstruction}. SRT yields significantly narrower wire target profiles than HMFA. The GEOFA and Transmit correction (see Remark \ref{rem:transmit_and_receive_correction}) further improves upon SRT. At $-5$\,dB, our correction significantly outperforms SRT above $15$\,mm depth and performs comparably below, closely matching the ground truth throughout. At $-15$\,dB, a difference between our method and SRT persists across all depths, and both deviate from the ground truth below $15$\,mm, although they still yield a substantial improvement over HMFA. The deviation from the ground truth observed for the wire target profile of target~$5$ and in the $-15$\,dB wire target width evaluation is attributed to RTB, which is used to obtain a higher resolution. In RTB, the data of several neighboring transmit events are combined coherently for each image pixel. Our correction is exact for the transmit event whose transmit line passes through the receive focus point. For the neighboring transmit events contributing to the same pixel, a residual transmit ToF error, smaller than the uncorrected HMFA error, remains. This residual error manifests as a broadening of the main lobe into the sidelobe. Correcting this effect would require incorporating the aberration correction already at the transmit stage or directly within the RTB algorithm, which would be a promising direction for future research.

Reverberations are visible as bright horizontal bands in the first $15$\,mm of the images in \cref{fig:9L_Fat_layer_images}. They are caused by multiple reflections between medium boundaries. In the wedge-aberrated setting, the primary reverberation arises from the wave reflecting at the silicone–phantom interface, traveling back to the transducer array inside the ultrasound transducer, and reflecting again toward the phantom. This artifact is most clearly visible in the GEOFA-corrected image due to the improved near-field image quality.

The results discussed in this section show that GEOFA combined with a transmit correction effectively corrects the aberrations introduced by a layer of differing sound speed. This correction achieved a significant image quality improvement compared to HMFA and the state-of-the-art aberration correction method of SRT. It should be noted, however, that a single experiment does not precisely establish which clinical situations require refraction-corrected focusing delays. Following the state-of-the-art literature \cite{jaeger2015computed,sanabria2018spatial,vraalstad2024coherence}, the SRT baseline contains no transmit correction, so the observed improvement is caused not only by the refraction-aware receive ToF calculation but also by the additional transmit correction of Remark~\ref{rem:transmit_and_receive_correction}. Moreover, the wedge aberrator represents a sound speed contrast of approximately $500\,\mathrm{m/s}$, which is comparable to that of a transducer cover or bone rather than to the milder contrast between soft tissues. Our simulations indicate that this is precisely the regime in which refraction matters. For the fat and tissue settings of \cref{sect_4_1_ToF_measurements}, the focusing delay errors of SRT remain below the significance threshold of Remark~\ref{rem:threshold}, whereas for the transducer cover and fetal imaging settings they exceed it by about one order of magnitude. A conclusive assessment of the practical relevance of refraction-aware techniques therefore requires further work across clinically realistic sound speed contrasts and geometries.

\section{Conclusion}\label{sect_5:Conclusion}

In this paper, we introduced a Geometrical Acoustics Focusing Algorithm (GEOFA) for ultrasound image quality improvement through a more precise ToF and focusing delay calculation. For this, we assumed a layered medium setting with known sound speeds in the medium layers and known, once continuously differentiable medium boundaries. Calculations based on geometrical acoustics were used to obtain the ToFs from sources to the desired focus point. Given these ToFs, the focusing delays, required to obtain the final image, could be calculated. Based on the underlying medium, we derived existence and uniqueness conditions for a solution of the method's underlying system of equations. Using the acoustic toolbox k-Wave, the ToFs and focusing delays for different layered media were simulated, and GEOFA was compared to SRT and the currently used homogeneous medium focusing algorithm. Numerical calculations were used to compare GEOFA with a similar algorithm derived in \cite{MINEO2022106747}. Finally, the ultrasound image improvements obtained by GEOFA in a phantom-based experimental setup were compared to SRT and the homogeneous medium focusing algorithm.

\section*{Acknowledgments and Support}

We would like to thank Dr. Stefan Denk, Dr. Peter Fosodeder and Dr. Martin Mienkina (GE HealthCare) for their valuable input, insightful feedback and crucial support in the experimental procedures. Additionally, we thank GE HealthCare for providing the equipment for our phantom-based experiments.

This research was funded in part by the Austrian Science Fund (FWF) SFB\\ 10.55776/F68 “Tomography Across the Scales”, project F6805-N36 (Tomography in Astronomy). For open access purposes, the authors have applied a CC BY public copyright license to any author-accepted manuscript version arising from this submission. The financial support by the Austrian Federal Ministry for Digital and Economic Affairs, the National Foundation for Research, Technology and Development and the Christian Doppler Research Association is gratefully acknowledged.

\appendix

\section{Existence and Uniqueness for the GEOFA System} \label{sect:6_Appendix}

This section discusses details of the existence and uniqueness theory from \cref{sect:3_3_MAFA_Existence_and_Uniqueness}.

\subsection{Existence of a Solution to the GEOFA System of Equations} \label{sect:6_1_Appendix_Existence_Result}

\begin{lemma}[No total reflection condition]\label{lemma:no_total_reflection_condition}
    For $1 \leq n \leq N-1$ and $1 \leq k \leq n$, let
    $P_0 = (x_0,z_0),$ and $P_k = (x_k,z_k)\in B_k$ be given. 
    Then, there occurs no total reflection on the $n$-th medium boundary in the point $P_n$ if and only if
    \begin{equation} \label{eq:MAFA_no_total_reflection_condition}
        \bigg| \frac{c_{n+1}}{c_n} \frac{ x_n - x_{n-1} + b_n'(x_n)(b_n(x_n) - b_{n-1}(x_{n-1}))}{\sqrt{1+b_n'(x_n)^2} \sqrt{(x_n - x_{n-1})^2 + (b_n(x_n) - b_{n-1}(x_{n-1}))^2}} \bigg| \leq 1 \, .
    \end{equation}
\end{lemma}

\begin{proof}
    We recall from \cref{sect:2_3_geometrical_acoustics} that no total reflection occurs in the point $P_n$ iff $|\frac{c_{n+1}}{c_n} \sin \theta_n| \leq 1$. Using the equations (\ref{eq:tangent_of_medium_boundary_equation}) and (\ref{eq:first_medium_boundary_eq}), the result follows. 
\end{proof}

\begin{figure}[ht!]
    \centering
    \includegraphics[width=.75\textwidth]{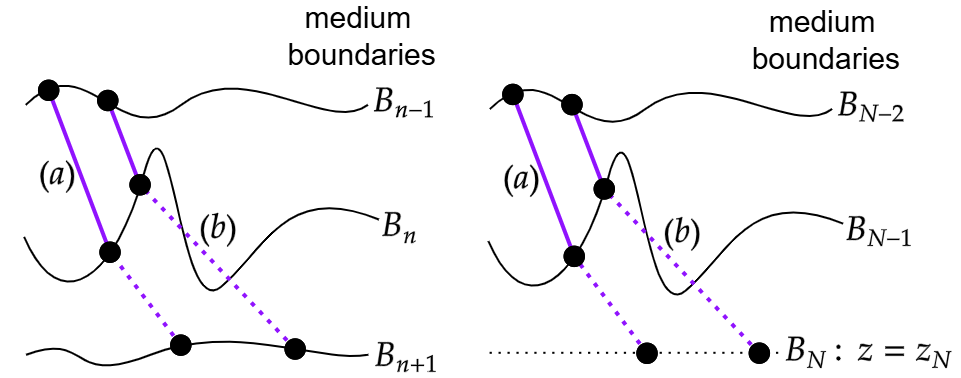}
    \caption{Setting in the proof of the unique intersection condition (Lemma~\ref{lemma:MAFA_unique_intersection_condition}). In both images, the acoustic ray $(a)$ has a unique intersection with the medium boundary, while the acoustic ray $(b)$ does not. Left: $n < N-1$, Right: $n = N-1$.}
    \label{fig:lemma_existence_of_a_mafa_path}
\end{figure}

Next, we derive a condition ensuring the unique intersection of an acoustic ray with each medium boundary.

\begin{lemma}[Unique intersection condition] \label{lemma:MAFA_unique_intersection_condition}
    Let $P_n \in B_n$ be the first intersection of the acoustic ray with $B_n$. Then the following statements hold:
     \begin{itemize}
         \item If the acoustic ray points straight down in the $(n+1)$-st medium layer, then there are no further intersections with $B_n$ between $P_n$ and $P_{n+1}$. 
        \item Otherwise, let $k_{n+1}$ be the slope of the acoustic ray in the $(n+1)$-st medium layer. Without considering further interactions with $B_n$, let $P_{n+1}$ be the first intersection of the acoustic ray with $B_{n+1}$ for $n < N-1$ or with the line $z = z_N$ for $n = N-1$, as visualized in \cref{fig:lemma_existence_of_a_mafa_path}. Then, there are no further intersections of the acoustic ray with $B_n$ between $P_n$ and $P_{n+1}$ if and only if 
        for all $x \in (\min\{ x_n,x_{n+1} \},\max\{ x_n,x_{n+1} \})$ there holds
        \begin{equation} \label{eq:MAFA_unique_medium_boundary_intersection_condition}
            b_n(x) - b_n(x_n) - k_{n+1} \cdot (x - x_n) = \int_{x_{n}}^x b_n'(y) - k_{n+1} \, dy \neq 0.
        \end{equation}  
        \item If \cref{eq:MAFA_unique_medium_boundary_intersection_condition} holds for all medium layers where the acoustic ray does not point straight down, then there are no additional intersections of the acoustic ray with the medium boundaries $B_1,\ldots,B_{N-1}$ other than $P_1,\ldots,P_{N-1}$.
     \end{itemize}
\end{lemma}

\begin{proof}
    If the acoustic ray points straight down in the $n+1$-st layer, no further intersections with the $n$-th medium boundary are possible due to the $z$-coordinate of the medium boundary being a function of the $x$-coordinate. 
    Otherwise, the acoustic ray does not intersect $B_n$ in between the points $P_n$ and $P_{n+1}$ if for all $x \in (\min\{ x_n,x_{n+1} \},$ $\max\{ x_n,x_{n+1} \})$ it holds that $ k_{n+1} \cdot (x - x_{n}) \neq b_n(x) - b_n(x_{n})$. This is equivalent to
        \[
            b_n(x) - b_n(x_{n}) - k_{n+1} \cdot (x - x_{n}) \neq 0\, ,
        \]
    which is the first part of condition (\ref{eq:MAFA_unique_medium_boundary_intersection_condition}). Using the fundamental theorem of calculus for $b_n \in C^1$, we obtain
    \[
        b_n(x) = b_n(x_{n}) + \int_{x_{n}}^x b_n'(y) \; dy \, ,
    \]
    and therefore the above term can equivalently be written in the form
    \[
        b_n(x) - b_n(x_n) - k_{n+1} (x - x_n) = \int_{x_{n}}^x b_n'(y) \, dy - k_{n+1} (x - x_n ) = \int_{x_{n}}^x b_n'(y) - k_{n+1} \, dy \, .
    \]
    Finally, since the segments between $P_n$ and $P_{n+1}$ partition the ray from $P_0$ to $P_N$ and no additional intersection exists within any such segment, the third assertion follows.
\end{proof}

At this point, we slightly change our perspective. Previously, we have fixed $P_0$ and $P_N$ and tried to find an acoustic ray through the medium, implicitly given by the GEOFA system of equations. Now, we fix $P_0$ and some point $\overline{P_1} = (\overline{x_1},b_1(\overline{x_1})) \in B_1$ and explicitly propagate the acoustic ray through the medium, which determines $\overline{P_n} \in B_n$ for $2 \leq n \leq N-1$. To emphasize this difference, these determined medium boundary points and their related variables are denoted with a bar. Furthermore, define $B_N:$ $z = z_N$ as the horizontal line containing the focus point. The intersection of the acoustic ray with $B_N$ is used for measuring the distance of the acoustic ray to the focus point. In the next lemma, we show that the intersections of the acoustic ray with the medium boundaries $B_2$, $\ldots$, $B_{N-1}$ and the line $B_N$ continuously depend on $\overline{x_1}$ for fixed $P_0$.

\begin{lemma} \label{lemma:medium_boundary_points_continuisly_dependence}
    Given a layered medium with medium boundaries $(B_n)_{n = 1}^{N-1}$ and a line $B_N:$ $z = z_N$. Let an acoustic ray originate in $P_0$ and pass through $(\overline{x_1},b_1(\overline{x_1})) = \overline{P_1} \in B_1$. Assume that for this acoustic ray, no total reflection and exactly one intersection with each medium boundary occurs. Then, all intersections $(\overline{P_n})_{n=1}^{N}$ of the acoustic ray with the medium boundaries $(B_n)_{n=1}^{N}$ and the line $B_N$ continuously depend on $\overline{x_1}$.
\end{lemma}
\begin{proof}

    We prove this lemma by using induction over the medium boundaries. The induction basis is trivial, as clearly $\overline{P_1} = (\overline{x_1},b_1(\overline{x_1}))$ continuously depends on $\overline{x_1}$. As induction hypothesis, we assume that $\overline{P_1}$, \ldots, $\overline{P_n} = (\overline{x_n},b_n(\overline{x_n}))$ continuously depend on $\overline{x_1}$. Therefore, also $\overline{x_1}$,\ldots, $\overline{x_n}$ continuously depend on $\overline{x_1}$. In the induction step, we show that $\overline{P_{n+1}} = (\overline{x_{n+1}},b_{n+1}(\overline{x_{n+1}}))$ continuously depends on $\overline{x_1}$ by directly calculating the acoustic ray's path through the $(n+1)$-st medium layer, as in the center and right panel of \cref{fig:exact_setting_in_MAFA_existence_proof}.
    
    Using the equations (\ref{eq:first_medium_boundary_eq}), (\ref{eq:tangent_of_medium_boundary_equation}) and Snell's law \cref{eq:Snells_Law} we obtain
    \[
        \overline{\theta_n'} = \arcsin \frac{c_n}{c_{n+1}} \frac{b_n'(\overline{x_n}) (b_n(\overline{x_n}) - b_{n-1}(\overline{x_{n-1}})) + (\overline{x_n} - \overline{x_{n-1}})}{\sqrt{b_n'(\overline{x_n})^2 + 1} \cdot \sqrt{(\overline{x_n} - \overline{x_{n-1}})^2 + (b_n(\overline{x_n}) - b_{n-1}(\overline{x_{n-1}}) )^2}} \, .
    \]
    The argument of the $\arcsin$ is in $[-1,1]$, since by assumption, no total reflection occurs for the acoustic ray, i.e., (\ref{eq:MAFA_no_total_reflection_condition}) is fulfilled. Therefore, $\overline{\theta_n'}$ continuously depends on $\overline{x_1}$. In the $(n+1)$-st layer, the acoustic ray in parameter form is given by 
    \begin{equation} \label{eq:appendix_equation_for_acoustic_ray}
        \overline{P_n} + t \cdot R(\overline{\theta'_n}) \overrightarrow{o_1},
    \end{equation}
    with the parameter $t \in \R^+$, $\overrightarrow{o_n} = (b_n'(\overline{x_n}),-1)^T$ being a vector pointing inside the $(n+1)$-st medium layer along the orthogonal of the tangent in $\overline{P_n}$, as visualized in \cref{fig:exact_setting_in_MAFA_existence_proof}, and $R(\overline{\theta'_n})$ being a rotation matrix:
    \[
        R(\overline{\theta'_n}) = \begin{pmatrix}
            \cos \overline{\theta'_{n}} & -\sin \overline{\theta'_{n}}\\
            \sin \overline{\theta'_{n}} & \cos \overline{\theta'_{n}}
        \end{pmatrix} \, .
    \]
    The acoustic ray (\ref{eq:appendix_equation_for_acoustic_ray}) is intersected with $B_{n+1}$ to obtain $\overline{P_{n+1}}$, which exists and is unique by assumption.
    
    Due to the form of the acoustic ray (\ref{eq:appendix_equation_for_acoustic_ray}), the continuity of $B_{n+1}$ and the induction hypothesis, $\overline{P_{n+1}}$ continuously depends on $\overline{x_1}$, which concludes our induction proof. 
\end{proof}

\subsection{Uniqueness of a Solution to the GEOFA System of Equations} \label{sect:6_2_Appendix_Uniqueness_Result}

This section provides a detailed proof of \cref{thm:uniqueness_GOAT_TFAE_statement} introduced in \cref{sect:3_3_MAFA_Existence_and_Uniqueness}.

\begin{definition}
    The \emph{change in Time of Flight} occurring if $P_n$ is moved along $B_n$ is defined as the partial derivative $\partial/\partial x_n \, \text{ToF}$ of the Time of Flight function, which is continuous as long as $P_n \neq P_{n-1}$ for all $n = 1,\ldots, N$.
\end{definition}

\begin{proof}[Proof of the equivalence of 1 and 3 in \cref{thm:uniqueness_GOAT_TFAE_statement}.]
    Since the $\text{ToF}$ function is continuously differentiable, its stationary points occur where its gradient vanishes. Fix $n$ that fulfills $1 \leq n \leq N-1$. Using $|P_n - P_{n-1}| = \sqrt{(\Delta x_n)^2 + (\Delta b_n(x_n))^2}$ and 
    \[
    \begin{aligned}
        \partial/\partial x_n &\text{ToF}(x_1,\ldots,x_{N-1})\\
        &=\frac{\Delta x_n + b_n'(x_n)  \Delta b_n(x_n)}{c_n \sqrt{(\Delta x_n)^2 + (\Delta b_n(x_n))^2}} - \frac{\Delta x_{n+1} + b_n'(x_n) \Delta b_{n+1}(x_{n+1})}{c_{n+1}\sqrt{(\Delta x_{n+1})^2 + (\Delta b_{n+1}(x_{n+1}))^2}} \, ,
    \end{aligned}
    \]
    the equation $\partial / \partial x_n \text{ToF}(x_1,\ldots,x_{N-1}) = 0$ can be rearranged to
    \[
        b_n'(x_n) \cdot \bigg( \frac{\Delta b_n(x_n)}{c_n | P_n - P_{n-1} |} - \frac{\Delta b_{n+1}(x_{n+1})}{c_{n+1} | P_{n+1} - P_n |} \bigg) = \frac{\Delta x_{n+1}}{c_{n+1} | P_{n+1} - P_n |} - \frac{\Delta x_n}{c_n | P_n - P_{n-1} |} \, .
    \]
    Multiplication with the common denominator yields
    \[
    \begin{aligned}
        b_n'(x_n) \cdot \big( c_{n+1} | P_{n+1} - P_n | \Delta b_n(x_n) &- c_n | P_n - P_{n-1} | \Delta b_{n+1}(x_{n+1})) \big)\\
        &= c_n | P_n - P_{n-1} | \Delta x_{n+1} - c_{n+1} | P_{n+1} - P_n | \Delta x_n \, .
    \end{aligned}
    \]
    Dividing by the factor on the left hand side now yields (\ref{eq:med_bound_slope_for_equal_ToF}).
\end{proof}
\begin{proof}[Proof of the equivalence of 2 and 3 in \cref{thm:uniqueness_GOAT_TFAE_statement}.]
    
    First, we show that Snell's law holds if (\ref{eq:med_bound_slope_for_equal_ToF}) is fulfilled. Inserting (\ref{eq:tangent_of_medium_boundary_equation}), (\ref{eq:first_medium_boundary_eq}) and (\ref{eq:second_medium_boundary_equation}) into Snell's law (\ref{eq:Snells_Law}) yields
    \[
    \begin{aligned}
        \frac{\Delta x_n + b_n'(x_n) \Delta b_n(x_n) }{c_n \sqrt{1 + b'(x_n)^2} \sqrt{(\Delta x_n)^2 + (\Delta b_n(x_n))^2}}\\
        &\hspace{-15mm}= \frac{\Delta x_{n+1} + b_n'(x_n) \Delta b_{n+1}(x_{n+1}) }{c_{n+1}\sqrt{1 + b'(x_n)^2} \sqrt{( \Delta  x_{n+1})^2 + (\Delta b_{n+1}(x_{n+1}))^2}}
    \end{aligned}
    \]
    By using the term $| P_{n+1} - P_n | = \sqrt{(\Delta x_{n+1})^2 + (\Delta b_{n+1}(x_{n+1}))^2}$ and multiplying with $\sqrt{1 + b'(x_n)^2}$, the denominators simplify, and we obtain
    \[
        \frac{\Delta x_n + b_n'(x_n) \Delta b_n(x_n) }{c_n  |P_{n-1} - P_{n}|} = \frac{\Delta x_{n+1} + b_n'(x_n) \Delta b_{n+1}(x_{n+1}) }{c_{n+1} |P_{n} - P_{n+1}|}
    \]
    Next, abbreviate $a_n := c_n | P_n - P_{n-1} |$ and $a_{n+1} := c_{n+1} | P_{n+1} - P_n |$ and multiply out:
    \[
        a_{n+1} (\Delta x_n + b_n'(x_n) \Delta b_n(x_n) ) = a_n (\Delta x_{n+1} + b_n'(x_n) \Delta b_{n+1}(x_{n+1})) \, .
    \]
    We assume that the slope of the $n$-th medium boundary, $b_n'(x_n)$, fulfills (\ref{eq:med_bound_slope_for_equal_ToF}). Inserting (\ref{eq:med_bound_slope_for_equal_ToF}) and multiplying by the denominator yields
    \[
    \begin{aligned}
        &a_{n+1} (\Delta x_n \cdot (a_{n+1} \Delta b_n(x_n) - a_n \Delta b_{n+1}(x_{n+1})) + \Delta b_n(x_n) \cdot (a_n \Delta x_{n+1} - a_{n+1} \Delta x_n)) \\
        &\qquad= a_n (\Delta x_{n+1} \cdot (a_{n+1} \Delta b_n(x_n) - a_n \Delta b_{n+1}(x_{n+1}))\\
        &\qquad \quad  + \Delta b_{n+1}(x_{n+1}) \cdot (a_n \Delta x_{n+1} - a_{n+1} \Delta x_n) \, ,
    \end{aligned}
    \]
    which finally simplifies to the following true statement
    \[
    \begin{aligned}
        a_n a_{n+1} ( \Delta b_n(x_n) \Delta x_{n+1} - \Delta &b_{n+1}(x_{n+1}) \Delta x_n)\\
        &= a_n a_{n+1} ( \Delta b_n(x_n) \Delta x_{n+1} - \Delta b_{n+1}(x_{n+1}) \Delta x_n)\, .
    \end{aligned}
    \]
    This shows that if $b_n'(x_n)$ fulfills (\ref{eq:med_bound_slope_for_equal_ToF}), Snell's law holds in $P_n$.

    For the second part of the proof, we show that Snell's law being fulfilled implies (\ref{eq:med_bound_slope_for_equal_ToF}). For this, we distinguish the cases $\theta_n = \theta_n'$ and $\theta_n \neq \theta_n'$. The case $\theta_n = \theta_n'$ is only possible if $\theta_n = \theta_n'= 0$, because Snell's law, $\theta_n \neq 0$ and $c_{n+1} \neq c_n$ imply $\theta_n \neq \theta_n'$. Therefore, the acoustic ray is orthogonal to the tangent in $P_n$ and no refraction occurs, as shown on the left of \cref{fig:appendix_6_2_image}. Since the two triangles with hypotenuses $|P_n - P_{n-1}|$ and $|P_{n+1} - P_{n-1}|$ are similar, there exists a constant $a \in \R^+$ s.t.
    \[
        |P_{n+1} - P_{n}| = a \, |P_n - P_{n-1}| \, , \quad \Delta x_{n+1} = a \; \Delta x_n \, , \quad \Delta b_{n+1}(x_{n+1}) = a \; \Delta b_n(x_n)\, .   
    \]
    Now, inserting into (\ref{eq:med_bound_slope_for_equal_ToF}) and using (\ref{eq:tangent_of_medium_boundary_equation}) yields 
    \begin{equation} \label{eq:appendix_6_2_reverse_direction_one_case}
        \tan \alpha = - \Delta x_{n+1} / \Delta b_{n+1}(x_{n+1}) \, .
    \end{equation}
    This concludes our proof of Snell's law implying (\ref{eq:med_bound_slope_for_equal_ToF}) in the case $\theta_n = \theta_n'$, since the fact that (\ref{eq:appendix_6_2_reverse_direction_one_case}) is correct can be obtained from the left panel of \cref{fig:appendix_6_2_image}. Note that one leg of the smaller triangle on the left panel of \cref{fig:appendix_6_2_image} is of length $-\Delta x_{n+1}$, because $x_{n+1} < x_n$ and therefore $\Delta x_{n+1} < 0$.
    
    \begin{figure}[ht!]
        \centering
        \includegraphics[width=.9\textwidth]{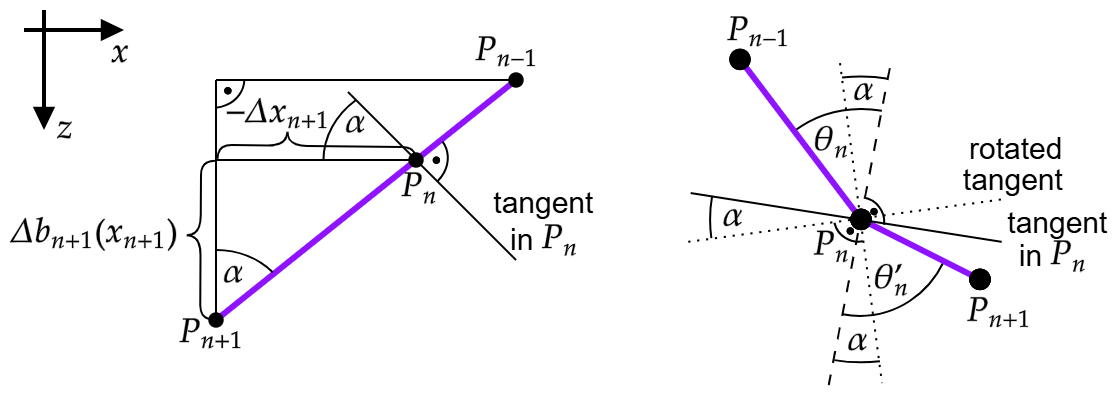}
        \caption{Proof of Snell's law implying (\ref{eq:med_bound_slope_for_equal_ToF}) for Left: $\theta_n = \theta_n'$, Right: $\theta_n \neq \theta_n'$}
        \label{fig:appendix_6_2_image}
    \end{figure}    

    In the case $\theta_n \neq \theta_n'$, we show that (\ref{eq:med_bound_slope_for_equal_ToF}) not being fulfilled implies that Snell's law does not hold. If (\ref{eq:med_bound_slope_for_equal_ToF}) is not fulfilled, the tangent in $P_n$ is rotated by some angle $\alpha$. Therefore, both $\theta_n$ and $\theta_n'$ are changed by the same angle $\alpha$, as visualized in \cref{fig:appendix_6_2_image}. To analyze this situation, we rewrite Snell's law to obtain
    \begin{equation} \label{eq:alternative_version_of_Snells_law}
        \sin \theta_n/\sin \theta_n'  = c_n/c_{n+1},
    \end{equation}
    and consider the function $f(\alpha) := \sin (\theta-\alpha)/\sin (\theta' - \alpha)$ for fixed angles $\theta, \theta'$ with $\theta'  - \alpha \neq k \cdot \frac{\pi}{2}$, $k \in \Z$. To finish this proof, we show that $f$ is strictly monotone, which means that (\ref{eq:alternative_version_of_Snells_law}) cannot be fulfilled for any value $\alpha \neq 0$, since we have shown in the first part of this proof that it holds for $\alpha = 0$.

    The strict monotonicity of $f$ can be shown using its derivative, which is given by
    \[
        f'(\alpha) = \big( \sin(\theta - \alpha) \cos(\theta' - \alpha) - \cos(\theta - \alpha) \sin(\theta' - \alpha) \big) /\sin^2(\theta' - \alpha) \, .
    \]
    With the angle addition theorem $\sin(x-y) = \sin(x) \cos(y) - \cos(x) \sin(y)$, we obtain
    \[
        f'(\alpha) = \sin(\theta - \theta') / \sin(\theta' - \alpha)^2 \neq 0
    \]
    for $\theta \neq \theta'$, implying that $f$ is strictly monotonous, as claimed.
\end{proof}

\section{k-Wave Simulation parameters used in \texorpdfstring{\cref{sect_4_1_ToF_measurements}}{Section 4.1}}
\label{sect:B_kWave_Simulation_Details}

The list of parameters used in the k-Wave simulations in \cref{sect_4_1_ToF_measurements} can be found in \cref{tab:simulation_parameters} and \cref{tab:geometry_simple}. Note that in k-Wave, the acoustic properties of the simulated medium have to be specified on a uniformly discretized spatial grid. For improved precision, we have linearly interpolated sound speeds of grid cells intersecting medium boundaries. Furthermore, the raw data obtained by the simulation consists of oscillating pressure curves. To obtain the envelope of the pressure curves, the absolute value of the Hilbert transform is used. The wave's arrival time is determined as the time point of the envelope's maximum. Finally, the precise emission time of the transmitted signal is evaluated in the same way, and subtracted from the arrival time to obtain the ToF.

\begin{table}[ht!]
\centering
\caption{Parameters used in the k-Wave simulations.}
\label{tab:simulation_parameters}
\begin{tabular}{llc}
\toprule
\multirow{6}{*}{Computational grid} & $x$-, $z$-discretization          & Settings 1-4: $5 \cdot 10^{-5}$\,m \\
                                    &                                   & Setting 5: $3 \cdot 10^{-5}$\,m           \\
                                    & \# of $x$-, $z$-grid points       & Settings 1-4: 729, 1715           \\
                                    &                                   & Setting 5: 1215, 3087             \\
                                    & Time discretization               & $2.5 \cdot 10^{-10}$\,s                    \\
                                    & Total simulation time             & $8 \cdot 10^{-5}$\,s                       \\
\midrule
Medium                              & \makecell{Sound speeds and\\medium geometries} & \multicolumn{1}{c}{\cref{tab:geometry_simple} and \cref{fig:ToF_simulations_Settings}} \\
\midrule
\multirow{3}{*}{Acoustic sources}   & Position                          & see \cref{fig:ToF_simulations_Settings} \\
                                    & \multirow{2}{*}{Emitted signal}   & \makecell{Gauss pulse,\\ \SI{5}{MHz} center frequency,} \\
                                    &                                   & fractional bandwidth $0.6$    \\
\midrule
\multirow{2}{*}{Sensors}            & Position                          & top of 1st layer                   \\
                                    & Recording                         & pressure                           \\
\bottomrule
\end{tabular}
\end{table}

\begin{table}[ht!]
\centering
\caption{Medium parameters for simulation Settings $1$ to $5$ of \cref{sect_4_1_ToF_measurements}. Thickness $d$ is measured at the center. Settings $1$ and $3$ use flat medium boundaries. The boundaries in Settings $2$, $4$, and $5$ are elliptic and the parameters $a$ and $b$ describe the elliptic medium boundary at the lower end of the respective layer.}
\label{tab:geometry_simple}
\begin{tabular}{cccccc}
\toprule
Setting& Layer & $c$ [\si{m/s}] & $d$ [\si{mm}] & $a$ [\si{mm}] & $b$ [\si{mm}] \\
\midrule
\multirow{2}{*}{1} & 1 & 1460 & 58 & \multicolumn{2}{c}{flat} \\
                   & 2 & 1540 & —  & —  & — \\[3pt]
\multirow{2}{*}{2} & 1 & 1460 & 60& 50 & 70 \\
                   & 2 & 1540 & —  & —  & — \\[3pt]
\multirow{2}{*}{3} & 1 & 1000 &  5 & \multicolumn{2}{c}{flat} \\
                   & 2 & 1540 & —  & —  & — \\[3pt]
\multirow{3}{*}{4} & 1 & 1540 &  10& 35 & 50 \\
                   & 2 & 2200 &  1 & 36 & 51 \\
                   & 3 & 1540 & —  & —  & — \\[3pt]
\multirow{8}{*}{5} & 1 & 1540 &  10& 35 & 50 \\
                   & 2 & 2200 &  1 & 36 & 51 \\
                   & 3 & 1600 &  2 & 38 & 53 \\
                   & 4 & 1460 &  20& 60 & 50 \\
                   & 5 & 1600 &  5 & 68 & 58 \\
                   & 6 & 1540 &  15& 80 & 60 \\
                   & 7 & 2200 &  2 & 78 & 58 \\
                   & 8 & 1540 & —  & —  & — \\
\bottomrule
\end{tabular}
\end{table}

\section{Sound Speed Reconstruction of the \cref{sect_4_3_fat_layer_image_improvement} Aberrator}
\label{sect:C_SoS_reconstruction}

\begin{figure}
    \centering
    \includegraphics[width=.66\linewidth]{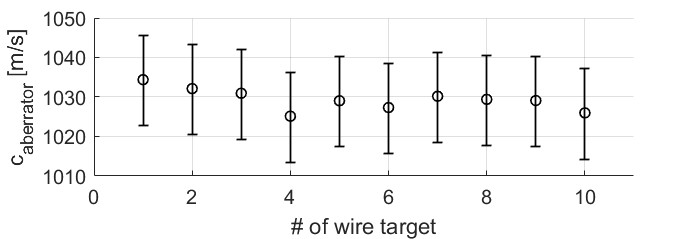}
    \caption{Reconstructed $c_{\operatorname{aberrator}}$ for eligible wire targets, ordered first from top to bottom and, for equal depth, left to right (see \cref{fig:9L_Fat_layer_images}). Error bars indicate uncertainty due to the 3D printer's layer thickness of $0.1\,mm$.}
    \label{fig:3D_printed_Aberrator_SoS_measurement}
\end{figure}

To reconstruct the sound speed $c_{\operatorname{aberrator}}$ of the aberrator used in \cref{sect_4_3_fat_layer_image_improvement}, the aberrated and ground truth images are compared. Inserting the cuboid aberrator between the transducer and the phantom shifts the wire targets in depth in the HMFA image, and this shift depends on $c_{\operatorname{aberrator}}$, the cuboid thickness, and the sound speed used in HMFA beamforming. 

For reconstructing $c_{\operatorname{aberrator}}$, all except two wire targets were eligible. The leftmost target in the horizontal row is excluded due to insufficient signal-to-noise ratio, while the rightmost target is excluded because of two closely spaced scatterers, which can also be observed in the wedge-aberrated setting in \cref{fig:9L_Fat_layer_images}. Both cases prevent reliable localization of the target center. The remaining targets yielded consistent reconstructions, suggesting that these exclusions do not bias the result, as evidenced by the low standard deviation. The reconstructed values of $c_{\text{aberrator}}$ are shown in \cref{fig:3D_printed_Aberrator_SoS_measurement}, yielding an arithmetic mean of
$c_{\text{aberrator}}= 1030.2 \pm 5.4 $\,m/s (standard deviation), with individual
values ranging from $1025.1$\,m/s to $1034.4$\,m/s. This value for $c_{\text{aberrator}}$ is within the range of sound speeds reported for silicone-rubber-based transducer lens materials~\cite{Yamashita2008}, making this experiment a realistic proxy for the aberrations they introduce.

\bibliographystyle{plain}
{\footnotesize
\bibliography{references}

\begin{thebibliography}{10}

\bibitem{Agarwal2016}
Mayur Agarwal, Arijit De, and Swapna Banerjee.
\newblock Architecture of a real-time delay calculator for digital beamforming
  in ultrasound system.
\newblock {\em IET Circuits Devices Syst.}, 10(4):322--329, 2016.

\bibitem{Alexandru2010}
Radu Alexandru.
\newblock Delay controller for ultrasound receive beamformer.
\newblock US Patent 7,804,736 B2, 2010.
\newblock Filed March 30, 2006, granted September 28, 2010.

\bibitem{Overview_of_Aberration_Correction_Methods}
R.~Ali, T.~Brevett, L.~Zhuang, H.~Bendjador, A.~S. Podkowa, S.~S. Hsieh,
  W.~Simson, S.~J. Sanabria, C.~D. Herickhoff, and J.~J. Dahl.
\newblock {Aberration correction in diagnostic ultrasound: A review of the
  prior field and current directions}.
\newblock {\em Z. Med. Phys.}, 33(3):267--291, 2023.

\bibitem{Sound_Speed_Estimation_using_SAFA_model}
R.~Ali, A.~V. Telichko, H.~Wang, U.~K. Sukumar, J.~G. Vilches-Moure,
  R.~Paulmurugan, and J.~J. Dahl.
\newblock {Local Sound Speed Estimation for Pulse-Echo Ultrasound in Layered
  Media}.
\newblock {\em IEEE Trans. Ultrason. Ferroelectr. Freq. Control},
  69(2):500--511, 2022.

\bibitem{Beaver1977}
William~L. Beaver.
\newblock Phase error effects in phased array beam steering.
\newblock In {\em Proceedings of the 1977 IEEE Ultrasonics Symposium}, pages
  264--267, 1977.

\bibitem{virtualSourceSAFT}
Nick Bottenus.
\newblock Comparison of virtual source synthetic aperture beamforming with an
  element-based model.
\newblock {\em J. Acoust. Soc. Am.}, 143(5):2801--2812, 2018.

\bibitem{Cass_Wu_Fink_closed_time_reversal_cavities}
D.~Cassereau, F.~Wu, and M.~Fink.
\newblock Limits of self-focusing using closed time-reversal cavities and
  mirrors-theory and experiment.
\newblock In {\em IEEE Symposium on Ultrasonics}, pages 1613--1618, 1990.

\bibitem{Ultrasound_propagation_velocity_important_for_image_Quality}
B.~Chauveau, C.~Auclair, A.~Legrand, R.~Mangione, L.~Gerbaud, F.~Vendittelli,
  L.~Boyer, and D.~L{\'e}mery.
\newblock {Improving image quality of mid-trimester fetal sonography in obese
  women: role of ultrasound propagation velocity}.
\newblock {\em Ultrasound Obstet. Gynecol.}, 52(6):769--775, 2018.

\bibitem{Chaves_nonimaging_optics}
J.~Chaves.
\newblock {\em {Introduction to Nonimaging Optics}}.
\newblock CRC Press, Boca Raton, 2nd edition, 2015.

\bibitem{Imaging_quality_deteriorates_for_obese_patients}
J.~S. Dashe, D.~D. McIntire, and D.~M. Twickler.
\newblock {Maternal obesity limits the ultrasound evaluation of fetal anatomy}.
\newblock {\em J. Ultrasound Med.}, 28(8):1025--1030, 2009.

\bibitem{Tissue_Properties_3}
Francis Duck.
\newblock {\em Physical Properties of Tissues: A comprehensive reference book}.
\newblock Academic Press, London, 1990.

\bibitem{FeynmanLectureSnell}
R.~Feynman, M.~A. Gottlieb, and R.~Pfeiffer.
\newblock Optics: The principle of least time, 1963.

\bibitem{Wave_Prop_Time_Rev_Randomly_Layered_Media}
J.-P. Fouque, J.~Garnier, G.~Papanicolaou, and K.~Solna.
\newblock {\em {Wave Propagation and Time Reversal in Randomly Layered Media}}.
\newblock Springer, New York, 2007.

\bibitem{Fritsch2006}
Carlos Fritsch, Montserrat Parrilla, Alberto Ib{\'a}{\~n}ez, Roberto~C.
  Giacchetta, and Oscar Mart{\'i}nez.
\newblock The progressive focusing correction technique for ultrasound
  beamforming.
\newblock {\em IEEE Trans. Ultrason. Ferroelectr. Freq. Control},
  53(10):1820--1831, 2006.

\bibitem{Tissue_Properties_1}
S.~A. Goss, R.~L. Johnston, and F.~Dunn.
\newblock {Comprehensive compilation of empirical ultrasonic properties of
  mammalian tissues}.
\newblock {\em J. Acoust. Soc. Am.}, 64(2):423--457, 1978.

\bibitem{Tissue_Properties_2}
S.~A. Goss, R.~L. Johnston, and F.~Dunn.
\newblock {Compilation of empirical ultrasonic properties of mammalian tissues.
  II}.
\newblock {\em J. Acoust. Soc. Am.}, 68(1):93--108, 1980.

\bibitem{Refraction_Based_SoS_estimation}
B.~H{\'e}riard-Dubreuil, A.~Besson, F.~Wintzenrieth, C.~Cohen-Bacrie, and J.-P.
  Thiran.
\newblock Refraction-based speed of sound estimation in layered media: An
  angular approach.
\newblock {\em IEEE Trans. Ultrason. Ferroelectr. Freq. Control},
  70(6):486--497, 2023.

\bibitem{ibrahim2017efficient}
Aya Ibrahim, Pascal~A Hager, Andrea Bartolini, Federico Angiolini, Marcel
  Arditi, Jean-Philippe Thiran, Luca Benini, and Giovanni De~Micheli.
\newblock Efficient sample delay calculation for 2-d and 3-d ultrasound
  imaging.
\newblock {\em IEEE Trans. Biomed. Circuits Syst.}, 11(4):815--831, 2017.

\bibitem{Jaeger_2022}
M.~Jaeger, P.~St{\"a}hli, N.~Korta~Martiartu, P.~Salemi~Yolgunlu, T.~Frappart,
  C.~Fraschini, and M.~Frenz.
\newblock Pulse-echo speed-of-sound imaging using convex probes.
\newblock {\em Phys. Med. Biol.}, 67(21):215016, 2022.

\bibitem{jaeger2015computed}
Michael Jaeger, Gerrit Held, Sara Peeters, Stefan Preisser, Michael Gr{\"u}nig,
  and Martin Frenz.
\newblock Computed ultrasound tomography in echo mode for imaging speed of
  sound using pulse-echo sonography: proof of principle.
\newblock {\em Ultrasound Med. Biol.}, 41(1):235--250, 2015.

\bibitem{Geometrical_Acoustics_Derivation}
J.~B. Keller.
\newblock {Geometrical Acoustics. I. The Theory of Weak Shock Waves}.
\newblock {\em J. Appl. Phys.}, 25(8):938--947, 1954.

\bibitem{HIFU_Focus_point_Displacement}
K.~Khoshgard, M.~Rezaei, and M.~Mousavi.
\newblock The correction of focal point displacement caused by the refraction
  of the beams in high-intensity focused ultrasound.
\newblock {\em J. Med. Signals Sens.}, 7:178--84, 2017.

\bibitem{Geometrical_Acoustics_for_solution_of_visualization_problems}
K.~G. Kvasnikov, A.~I. Soldatov, I.~O. Bolotina, Kh.~M. Krening, and A.~A.
  Potapenko.
\newblock The use of geometrical acoustics for the solution of visualization
  problems.
\newblock {\em Russ. J. Nondestruct. Test.}, 49(11):625--630, 2013.

\bibitem{liu1997propagation}
D.-L. Liu and R.~C. Waag.
\newblock Propagation and backpropagation for ultrasonic wavefront design.
\newblock {\em IEEE Trans. Ultrason. Ferroelectr. Freq. Control}, 44(1):1--13,
  1997.

\bibitem{4_5_per_cent_Salt_water_at_room_temperature_1540_Soundspeed}
K.~V. Mackenzie.
\newblock Nine-term equation for sound in the oceans.
\newblock {\em J. Acoust. Soc. Am.}, 70(3):807--812, 1981.

\bibitem{Magnin1981}
Paul~A. Magnin, Olaf~T. von Ramm, and Frederick~L. Thurstone.
\newblock Delay quantization error in phased array images.
\newblock {\em IEEE Trans. Sonics Ultrason.}, SU-28(5):305--310, 1981.

\bibitem{mast1999simulation}
T.~D. Mast, L.~M. Hinkelman, L.~A. Metlay, M.~J. Orr, and R.~C. Waag.
\newblock Simulation of ultrasonic pulse propagation, distortion, and
  attenuation in the human chest wall.
\newblock {\em J. Acoust. Soc. Am.}, 106(6):3665--3677, 1999.

\bibitem{mast1997simulation}
T.~D. Mast, L.~M. Hinkelman, M.~J. Orr, V.~W. Sparrow, and R.~C. Waag.
\newblock Simulation of ultrasonic pulse propagation through the abdominal
  wall.
\newblock {\em J. Acoust. Soc. Am.}, 102(2):1177--1190, 1997.

\bibitem{ToF_Distortions_caused_by_Large_Scale_Tissue_layers}
T.~D. Mast, L.~M. Hinkelman, M.~J. Orr, and R.~C. Waag.
\newblock The effect of abdominal wall morphology on ultrasonic pulse
  distortion. part ii. simulations.
\newblock {\em J. Acoust. Soc. Am.}, 104:3651--3664, 1998.

\bibitem{MINEO2022106747}
C.~Mineo, D.~Cerniglia, and E.~Mohseni.
\newblock {Solving ultrasonic ray tracing in parts with multiple material
  layers through Root-Finding methods}.
\newblock {\em Ultrasonics}, 124:106747, 2022.

\bibitem{MINEO2021106330}
Carmelo Mineo, David Lines, and Donatella Cerniglia.
\newblock Generalised bisection method for optimum ultrasonic ray tracing and
  focusing in multi-layered structures.
\newblock {\em Ultrasonics}, 111:106330, 2021.

\bibitem{Odegaard1995}
L.~A. {\O}degaard.
\newblock {\em Phase aberration correction in medical ultrasound imaging}.
\newblock {PhD} thesis, Norwegian Institute of Technology, University of
  Trondheim, Trondheim, Norway, 1995.

\bibitem{So_you_think_you_can_DAS}
V.~Perrot, M.~Polichetti, F.~Varray, and D.~Garcia.
\newblock {So you think you can DAS? A viewpoint on delay-and-sum beamforming}.
\newblock {\em Ultrasonics}, 111:106309, 2021.

\bibitem{Focusing_error_source}
D.~K. Peterson and Gordon~S. Kino.
\newblock Real-time digital image reconstruction: A description of imaging
  hardware and an analysis of quantization errors.
\newblock {\em IEEE Trans. Sonics Ultrason.}, 31(4):337--351, 1984.

\bibitem{pinton2011sources}
G.~F. Pinton, G.~E. Trahey, and J.~J. Dahl.
\newblock Sources of image degradation in fundamental and harmonic ultrasound
  imaging using nonlinear, full-wave simulations.
\newblock {\em IEEE Trans. Ultrason. Ferroelectr. Freq. Control},
  58(4):754--765, 2011.

\bibitem{Acoustic_Impedance_Matching}
Vivek~T. Rathod.
\newblock A review of acoustic impedance matching techniques for piezoelectric
  sensors and transducers.
\newblock {\em Sensors}, 20(14), 2020.

\bibitem{USTB}
A.~Rodriguez-Molares, O.~M.~H. Rindal, O.~Bernard, A.~Nair, M.~A. Lediju~Bell,
  H.~Liebgott, A.~Austeng, and L.~Lovstakken.
\newblock The ultrasound toolbox.
\newblock In {\em 2017 IEEE International Ultrasonics Symposium (IUS)}, pages
  1--4, 2017.

\bibitem{sanabria2018spatial}
Sergio~J Sanabria, Ece Ozkan, Marga Rominger, and Orcun Goksel.
\newblock Spatial domain reconstruction for imaging speed-of-sound with
  pulse-echo ultrasound: simulation and in vivo study.
\newblock {\em Phys. Med. Biol.}, 63(21):215015, 2018.

\bibitem{geometrical_room_acoustic_modeling}
Lauri Savioja and U.~Peter Svensson.
\newblock Overview of geometrical room acoustic modeling techniques.
\newblock {\em J. Acoust. Soc. Am.}, 138:708--730, 2015.

\bibitem{shen2010computational}
Y.-T. Shen, M.~I. Daoud, and J.~C. Lacefield.
\newblock Computational models of distributed aberration in ultrasound breast
  imaging.
\newblock {\em IEEE Trans. Ultrason. Ferroelectr. Freq. Control},
  57(12):2627--2636, 2010.

\bibitem{Multilayer_Synthetic_Aperature_Focusing}
M.~Skjelvareid and Y.~Birkelund.
\newblock Ultrasound imaging using multilayer synthetic aperture focusing.
\newblock In {\em Proceedings of the ASME 2010 Pressure Vessels and Piping
  Division/K-PVP Conference}, pages 379--387, 2010.

\bibitem{Steinberg1992}
Bernard~D. Steinberg.
\newblock Digital beamforming in ultrasound.
\newblock {\em IEEE Trans. Ultrason. Ferroelectr. Freq. Control},
  39(6):716--721, 1992.

\bibitem{szabo2013diagnostic}
T.~L. Szabo.
\newblock {\em Diagnostic Ultrasound Imaging: Inside Out}.
\newblock Elsevier Academic Press, Burlington, MA, 2013.

\bibitem{Thomenius_Evolution_Beamforming}
K.~E. Thomenius.
\newblock Evolution of ultrasound beamformers.
\newblock In {\em 1996 IEEE Ultrasonics Symposium}, volume~2, pages 1615--1622,
  1996.

\bibitem{First_kWave_Paper}
B.~Treeby and B.~T. Cox.
\newblock {k-Wave: MATLAB toolbox for the simulation and reconstruction of
  photoacoustic wave fields}.
\newblock {\em J. Biomed. Opt.}, 15:021314, 2010.

\bibitem{villarino2007decartesperfectlens}
M.~B. Villarino.
\newblock Decartes' perfect lens, 2007.

\bibitem{vraalstad2024coherence}
Anders~Emil Vr{\aa}lstad, Peter Fosodeder, Karin~Ulrike Deibele, Siri~Ann
  Nyrnes, Ole Marius~Hoel Rindal, Vibeke Skoura-Torvik, Martin Mienkina, and
  Svein-Erik M{\aa}s{\o}y.
\newblock {Coherence-Based Sound Speed Aberration Correction---With Clinical
  Validation in Fetal Ultrasound}.
\newblock {\em IEEE Trans. Med. Imaging}, 45(7):4047--4057, 2026.

\bibitem{TD_Cover_reconstruction_and_influence_on_image_quality}
Rick Waasdorp, David Maresca, and Guillaume Renaud.
\newblock Assessing transducer parameters for accurate medium sound speed
  estimation and image reconstruction.
\newblock {\em IEEE Trans. Ultrason. Ferroelectr. Freq. Control},
  71(10):1233--1243, 2024.

\bibitem{Overview_Ultrasound_Imaging}
P.~N.~T. Wells.
\newblock Ultrasound imaging.
\newblock {\em J. Biomed. Eng.}, 10(6):548--554, 1988.

\bibitem{SoS_estimate_steering}
Di~Xiao, Pat De~la Torre, and Alfred C.~H. Yu.
\newblock Real-time speed-of-sound estimation in vivo via steered plane wave
  ultrasound.
\newblock {\em IEEE Trans. Ultrason. Ferroelectr. Freq. Control},
  71(6):673--686, 2024.

\bibitem{Yamashita2008}
Y~Yamashita, Y~Hosono, and K~Itsumi.
\newblock Low-attenuation acoustic silicone lens for medical ultrasonic array
  probes.
\newblock In {\em Piezoelectric and acoustic materials for transducer
  applications}, pages 161--177. Springer, 2008.

\end{thebibliography}
}

\end{document}